\documentclass[graybox]{svmult}

\usepackage{color}
\usepackage{graphicx}
\usepackage{amsmath}
\usepackage{amssymb}
\usepackage{amsfonts}
\usepackage{array}
\usepackage{multirow}
\usepackage[colorlinks,bookmarksopen,bookmarksnumbered,citecolor=red,urlcolor=red]{hyperref}
\newcolumntype{C}[1]{>{\centering\let\newline\\\arraybackslash\hspace{0pt}}m{#1}}
\DeclareMathOperator*{\spann}{span}
\DeclareMathOperator*{\spec}{spec}
\usepackage{enumitem}
\usepackage{setspace}

\usepackage{mathptmx}
\usepackage{makeidx}
\usepackage{multicol} 
\usepackage[bottom]{footmisc}
\makeindex 

\begin{document}

\title{On the optimality of shifted Laplacian in a class of polynomial preconditioners for the Helmholtz equation}
\titlerunning{On a class of polynomial preconditioners for the Helmholtz equation}
\author{Siegfried Cools, Wim Vanroose}
\institute{Siegfried Cools \at Department of Mathematics and Computer Science, University of Antwerp, Middelheimlaan 1, 2020 Antwerp, Belgium, \email{siegfried.cools@uantwerp.be}
\and Wim Vanroose \at Department of Mathematics and Computer Science, University of Antwerp, Middelheimlaan 1, 2020 Antwerp, Belgium \email{wim.vanroose@uantwerp.be}}

\abstract{ 
This paper introduces and explores a class of polynomial preconditioners for the Helmholtz equation, denoted as expansion preconditioners $EX(m)$, that form a direct generalization to the classical complex shifted Laplace (CSL) preconditioner. The construction of the $EX(m)$ preconditioner is based on a truncated Taylor series expansion of the original Helmholtz operator inverse. The expansion preconditioner is shown to significantly improve Krylov solver convergence rates for growing values of the number of series terms $m$. However, the addition of multiple terms in the expansion also increases the computational cost of applying the preconditioner. A thorough cost-benefit analysis of the addition of extra terms in the $EX(m)$ preconditioner proves that the CSL or $EX(1)$ preconditioner is the most efficient member of the expansion preconditioner class for general practical and solver problem settings. Additionally, possible extensions to the expansion preconditioner class that further increase preconditioner efficiency are suggested, and numerical experiments in 1D and 2D are presented to validate the theoretical results.
}


\maketitle

\textbf{Summary.}  
This paper introduces and explores a class of polynomial preconditioners for the Helmholtz equation, denoted as expansion preconditioners $EX(m)$, that form a direct generalization to the classical complex shifted Laplace (CSL) preconditioner. The construction of the $EX(m)$ preconditioner is based on a truncated Taylor series expansion of the original Helmholtz operator inverse. The expansion preconditioner is shown to significantly improve Krylov solver convergence rates for growing values of the number of series terms $m$. However, the addition of multiple terms in the expansion also increases the computational cost of applying the preconditioner. A thorough cost-benefit analysis of the addition of extra terms in the $EX(m)$ preconditioner proves that the CSL or $EX(1)$ preconditioner is the most efficient member of the expansion preconditioner class for general practical problem and solver settings. Additionally, possible extensions to the expansion preconditioner class that further increase preconditioner efficiency are suggested, and numerical experiments in 1D and 2D are presented to validate the theoretical results.

\section{Introduction}
\label{sec:introduction}

\subsection{Overview of recent developments}

The propagation of waves through any material is often mathematically modeled by the Helmholtz equation, which represents the time-independent waveforms in the frequency domain. For high wavenumbers, i.e. high spatial frequencies, the sparse linear system that results from the discretization of this PDE is distinctly indefinite, causing most of the classic direct and iterative solution methods to perform poorly. Over the past few years, many different Helmholtz solution methods have been proposed, an overview of which can be found in \cite{ernst2012difficult}. Krylov subspace methods like GMRES \cite{saad1986gmres} or BiCGStab \cite{van1992bicg} are known for their robustness and are hence frequently used for the solution of Helmholtz problems \cite{lairdpreconditioned,bayliss1983iterative,simoncini2007recent,osei2010preconditioning}. However, due to the indefinite nature of the problem, Krylov methods are generally not efficient as Helmholtz solvers without the inclusion of a suitable preconditioner. 

\begin{figure}[t!] 
\begin{center}
\includegraphics[width=0.70\textwidth]{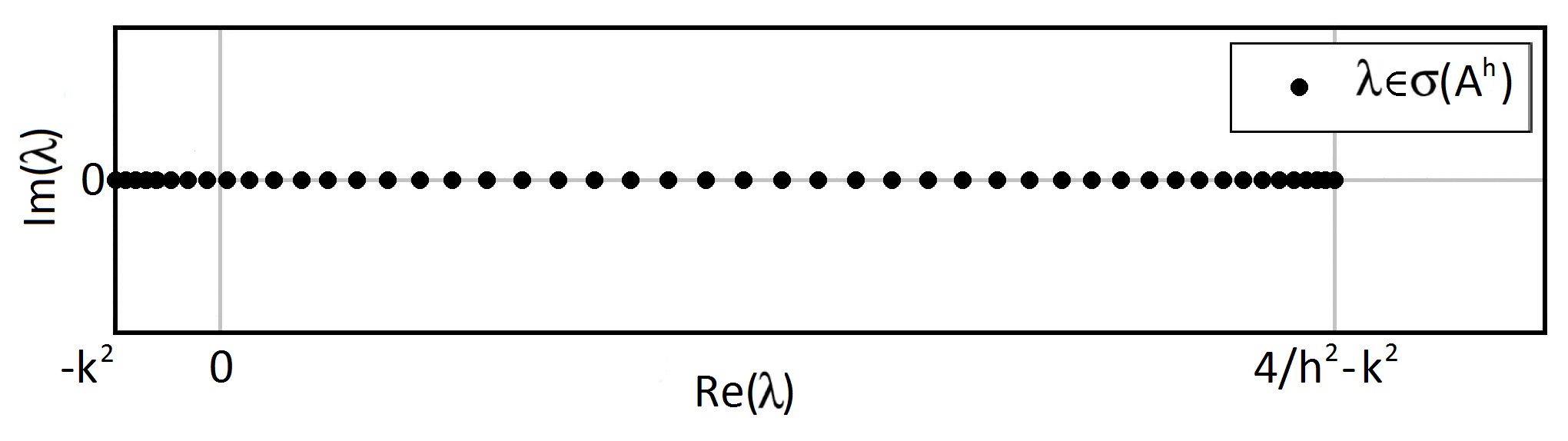}
\includegraphics[width=0.70\textwidth]{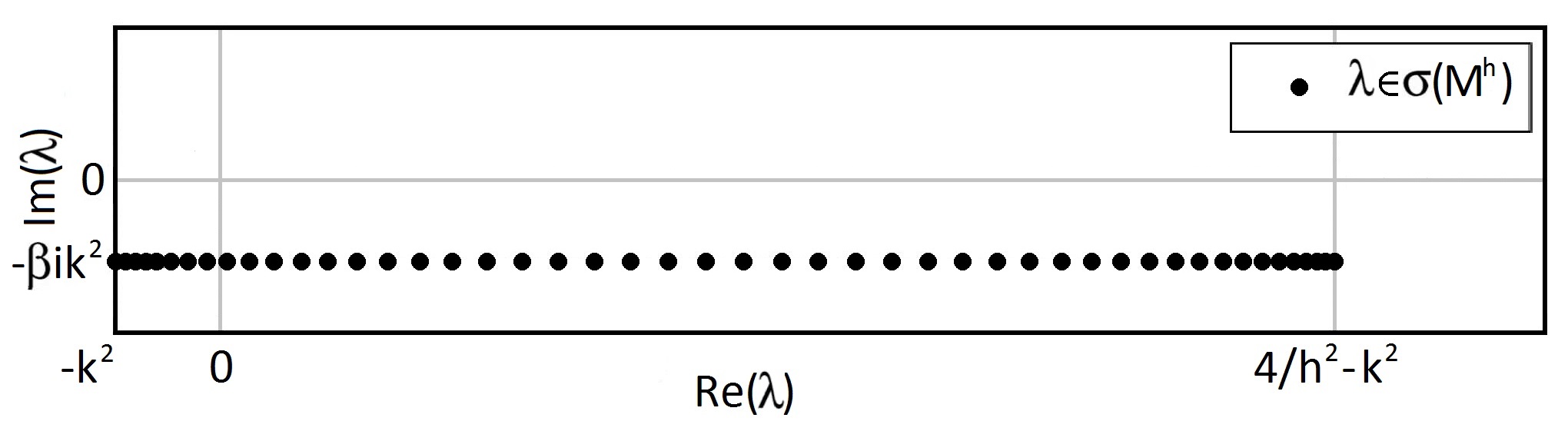}
\caption{Top: spectrum of the 1D Helmholtz operator discretized using second order finite differences on a $48 \times 48$ equidistant grid with standard homogeneous Dirichlet boundary conditions. Bottom: spectrum of the corresponding Complex Shifted Laplacian (CSL) operator.}
\label{fig:spectra_HH_and_CSL}
\end{center}
\end{figure}

In this chapter we focus on the class of so-called shifted Laplace preconditioners, which were introduced in \cite{made2001incomplete} and \cite{erlangga2004class}, and further analyzed in \cite{erlangga2006novel,erlangga2006comparison}. It was shown in the literature that contrary to the original discretized Helmholtz system, the complex shifted Laplace (CSL) system (or damped Helmholtz equation) can be solved efficiently using iterative methods \cite{reps2010indefinite,ernst2012difficult}. Originally introduced by Fedorenko in \cite{fedorenko1964speed}, multigrid methods \cite{brandt1977multi,brandt1986multigrid,stüben1982multigrid,briggs2000multigrid,trottenberg2001multigrid,elman2002multigrid} have been proposed as scalable solution methods for the shifted Laplace system in the literature \cite{erlangga2006comparison}. Typically only one multigrid V-cycle on the CSL system yields a sufficiently good approximate inverse, which can then be used as a preconditioner to the original Helmholtz system \cite{erlangga2006novel,erlangga2008multilevel,reps2012analyzing,sheikh2013convergence}. 
The main concept behind the shifted Laplace preconditioner is deceivingly simple: by shifting the spectrum of the Helmholtz operator down into the complex plane, close-to-zero eigenvalues (leading to near-singularity) that possibly destroy the iterative solver convergence are avoided, as illustrated in Figure \ref{fig:spectra_HH_and_CSL}. Nevertheless, the results in this work will show that this apparent simplicity is exactly what makes the CSL preconditioner into a powerful tool for the iterative solution of the Helmholtz equation.

\subsection{Outline of this chapter}

This study presents a generalization of the class of shifted Laplace preconditioners, which is based on a Taylor series expansion \cite{morse1953methods,bhatia1997matrix,arfken2011mathematical} of the original Helmholtz operator inverse around a complex shifted Laplacian operator. This formulation relates the original Helmholtz inverse to an infinite sum of shifted Laplace problems. By truncating the series we are able to define a class of so-called expansion preconditioners, denoted by $EX(m)$, where the number of terms $m$ in the expansion is a parameter of the method. The expansion preconditioner directly generalizes the classic complex shifted Laplace preconditioner, since the CSL preconditioner appears as the operator $EX(1)$, i.e.\ the first term in the Taylor expansion. 

Using a spectral analysis \cite{erlangga2006novel,vangijzen2007spectral,cools2013local}, the incorporation of additional series terms in the $EX(m)$ preconditioning polynomial is shown to greatly improve the spectral properties of the preconditioned system. When used as a preconditioner the $EX(m)$ operator hence allows for a significant reduction of the number of outer Krylov steps required to solve the Helmholtz problem for growing values of $m$. However, the addition of multiple terms in the expansion also increases the computational cost of applying the preconditioner, since each additional series term comes at the cost of one extra shifted Laplace operator inversion. The performance trade-off between the reduction of the number of outer Krylov iterations and the cost of additional terms (CSL inversions) in the preconditioner polynomial is analyzed in-depth. Furthermore, several theoretical extensions to the expansion preconditioner are introduced to improve preconditioner efficiency. These extensions provide the reader with supplementary insights into Helmholtz preconditioning. The proposed methods show similarities to the research on flexible Krylov methods \cite{saad1993fgmres,simoncini2003flexible} and more specifically the work on multi-preconditioned GMRES \cite{greif2011multi}.

A variety of numerical experiments are performed to validate the $EX(m)$ preconditioner and illustrate the influence of the number of series terms $m$ on convergence. Performance and computational cost of CSL- and $EX(m)$-preconditioned BiCGStab \cite{van1992bicg} are compared for one- and two-dimensional Helmholtz model problems with absorbing boundary conditions. These absorbing boundaries are implemented using Exterior Complex Scaling (ECS) \cite{aguilar1971class,simon1979definition,reps2012analyzing}, a technique which has been related to Perfectly Matched Layers (PMLs) \cite{berenger1994perfectly} by Chew and Weedon \cite{chew19943d}.

The remainder of this chapter is organized as follows. In Section \ref{sec:framework} we outline the theoretical framework for this work and we introduce the expansion preconditioner class $EX(m)$. Following its formal definition, an overview of the theoretical, numerical and computational properties of the expansion preconditioner is given. Section \ref{sec:extensions} presents several possible extensions to the proposed $EX(m)$ preconditioner, which are shown to improve preconditioner efficiency even further. The new preconditioner class is validated in Section \ref{sec:experiments}, where it is applied to a 1D and 2D Helmholtz benchmark problem. A spectral analysis confirms the asymptotic exactness of the expansion preconditioner as the number of terms $m$ grows towards infinity. Additionally, experiments are performed to compare the efficiency and computational cost of the $EX(m)$ preconditioner for various values of $m$. Conclusions and a short discussion on the results are formulated in Section \ref{sec:conclusions}.

\section{The expansion preconditioner} 
\label{sec:framework}

In this section we introduce the general framework for the construction of the expansion preconditioner. Starting from the notion of the classic complex shifted Laplace operator, we define the class of expansion preconditioners based on a Taylor expansion of the original Helmholtz operator inverse around a shifted Laplace problem with an arbitrary shift parameter. The definition of the expansion preconditioner is followed by an overview of the fundamental analytical and computational properties of the new preconditioner class. 

\subsection{The complex shifted Laplacian preconditioner}

In this work we aim to construct an efficient solution method for the $d$-dimensional Helmholtz equation 
\begin{equation} \label{eq:hh_sys}
( -\Delta - k^2(\bold{x}) ) \, u(\bold{x}) = f(\bold{x}), \quad \bold{x} \in \Omega \subset \mathbb{R}^d,
\end{equation}
with outgoing wave boundary conditions
\begin{equation}
u = \textrm{outgoing on } \partial \Omega,
\end{equation}
where $-k(\bold{x})^2$ is a distinctly negative shift. Here $k \in \mathbb{R}$ designates the wavenumber, which will be assumed to be a spatially independent constant throughout most of this work for simplicity. However, note that the definitions in this section do not depend on this assumption. The above equation is discretized using a finite difference, finite element or finite volume scheme of choice, yielding a system of linear equations of the general form
\begin{equation}\label{eq:hh_sys_dis}
Au = f,
\end{equation}
where the matrix operator $A$ represents a discretization of the Helmholtz operator $A \overset{d}{=} (-\Delta -k^2)$. It has been shown in the literature that iterative methods in general, and multigrid in particular, fail at efficiently solving the discretized Helmholtz system \eqref{eq:hh_sys_dis} due to the indefiniteness of the operator $A$ \cite{elman2002multigrid,ernst2012difficult}. However, the addition of a complex shift in the Helmholtz system induces a damping. This allows for a more efficient solution of the resulting system, which is known as the complex shifted Laplacian (CSL)
\begin{equation} \label{eq:csl_sys}
( -\Delta - (1+\beta i)k^2(\bold{x}) ) \, u(\bold{x}) = f(\bold{x}), \quad \bold{x} \in \Omega \subset \mathbb{R}^d,
\end{equation}
where $\beta \in \mathbb{R}^{+}$ is the complex shift (or damping) parameter that is conventionally chosen to be positive \cite{erlangga2004class,erlangga2006comparison}. The discretized shifted Laplace system is denoted by
\begin{equation} \label{eq:csl_sys_dis}
M u = f,
\end{equation}
where $M$ is the discretization of the complex shifted Helmholtz operator $M \overset{d}{=}(-\Delta - (1+\beta i)k^2)$. 

It is well-known that this system can be solved using multigrid when the shift parameter $\beta$ is sufficiently large \cite{erlangga2004class,erlangga2006novel,erlangga2008multilevel}. Furthermore, it has been shown in the literature that the complex shift parameter $\beta$ should be chosen at least as large as the critical value $\beta_{\min}$. For a multigrid V-cycle with the standard linear interpolation and full weighting restriction operators and a traditional $\omega$-Jacobi or Gauss-Seidel type smoothing scheme, the rule-of-thumb value for $\beta_{\min}$ was shown to be $0.5$ for a V(1,0)-cycle \cite{erlangga2006novel}, and lies roughly around $0.6$ when solving the CSL system using a V(1,1) multigrid cycle \cite{cools2013local}. Note that the latter value for the shift will be commonly used throughout this chapter.

The exact solution to \eqref{eq:csl_sys} is a damped waveform that is fundamentally different from the solution to the original Helmholtz system \eqref{eq:hh_sys}. However, the inverse of the shifted matrix operator $M$ can be used as a preconditioner to the original system. This preconditioning technique is known as the complex shifted Laplace preconditioner, and was shown to perform well as a preconditioner for Helmholtz problems, see \cite{erlangga2006novel,erlangga2008multilevel,vangijzen2007spectral}.

\subsection{The class of expansion preconditioners}

\subsubsection{General criteria for preconditioner efficiency} \label{sec:measure}

The aim of this work is to extend the existing class of shifted Laplace preconditioners to obtain a more efficient preconditioner to the original Helmholtz system \eqref{eq:hh_sys_dis}. Moreover, we would like to construct a preconditioner $M$ such that $M \approx A$ or, as an equivalent measure, we require that the eigenvalues of the preconditioned operator $M^{-1}A$ are concentrated around one, i.e.
\begin{equation} \label{eq:effec2}
\spec(M^{-1} A) \approx 1.
\end{equation}
In the context of an efficient iterative solution, the condition number $\kappa = \kappa(M^{-1} A)$ of the preconditioned system is often related to the number of Krylov iterations \cite{liesen2012krylov}. Although this relation is somewhat heuristic in the context of non-normal matrices \cite{trefethen2005spectra}, we believe that it provides an insightful intuition on preconditioner efficiency. The above requirement \eqref{eq:effec2} can hence be broadened by requiring that the condition number $\kappa$ approximately equals one, i.e.
\begin{equation} \label{eq:effec}
\kappa(M^{-1} A) \approx 1.
\end{equation}
We trivially note that by the above characterizations the best preconditioner to the Helmholtz system \eqref{eq:hh_sys_dis} is (a good approximation to) the original operator $A$ itself. However, the discrete operator $A$ is generally close to singular and hence cannot be easily inverted in practice. On the other hand, given a sufficiently large complex shift, the CSL system \eqref{eq:csl_sys} can be approximately inverted using e.g.\ a (series of) multigrid V-cycle(s). Note that the complex shifted Laplacian is generally not a very precise approximation to the original Helmholtz operator unless the shift parameter $\beta$ is very small, which in practice is not achievable due to stability requirements on the preconditioner inversion, see \cite{cools2013local}.

Following the general idea of preconditioning, the optimal Helmholtz preconditioner $M_{\text{opt}}$ thus ideally satisfies the following two key properties, which are inspired by analogous conditions that were formulated in \cite{gander2015applying}:
\vspace{0.1cm} \\ 
\fbox{ \begin{minipage}{0.96\textwidth}	{\setstretch{1.2}
	\begin{itemize}[leftmargin=0.8cm]
	\item[(P1)] $M_{\text{opt}}^{-1}$ is a good approximation to the exact inverse of the original 
							Helmholtz operator $A^{-1}$, such that condition \eqref{eq:effec2} is satisfied, 
							i.e.~the spectrum of the preconditioned operator $M_{\text{opt}}^{-1} A$ is clustered around one, 
	\end{itemize} }
\end{minipage}
}
\vspace{0.1cm} \\
\fbox{ \begin{minipage}{0.96\textwidth} {\setstretch{1.2}
	\begin{itemize}[leftmargin=0.8cm]
	\item[(P2)] for any given vector $v$, $M_{\text{opt}}^{-1}v$ can be efficiently computed iteratively. 
							This implies a good approximation to $M_{\text{opt}}^{-1}v$ is found after a `moderate' 
							number of iterations of the chosen method, with a manageable cost per iteration. 
	\end{itemize} }
\end{minipage}
}
\vspace{0.1cm} \\
In the context of this chapter condition (P2) is satisfied if $M_{\text{opt}}^{-1}$ can be formulated in terms of shifted Laplace operator inverses with a shift parameter $\beta$ that is sufficiently large to ensure a stable iterative solution of the shifted Laplace inverses. Indeed, given that the shift parameter $\beta$ is sufficiently large, a good approximation to the CSL inverse can be computed using e.g.~only one multigrid V-cycle, see \cite{erlangga2006novel}.

Note that due to the strong indefiniteness of the Helmholtz operator conditions (P1) and (P2) are generally incompatible. The classic CSL preconditioner trivially satisfies the second condition given that the shift parameter $\beta$ is large enough, however, the first condition is typically violated when $\beta$ is large. In the following we aim at constructing a preconditioning scheme based upon the shifted Laplace preconditioner, which effectively satisfies \emph{both} of the above conditions. 

\subsubsection{Taylor series expansion of the inverse Helmholtz operator}

The complex shifted Laplace preconditioning operator $M$ can be written more generally as
\begin{equation} \label{eq:hhcharac}
M(\beta) \overset{d}{=} -\Delta - k^2(\bold{x}) + P(\beta,\bold{x}),
\end{equation}
where $P(\beta,\bold{x})$ is a possibly spatially dependent linear operator in the shift parameter $\beta \in \mathbb{R}^{+}$, satisfying $P(0,\bold{x}) = 0$, such that $M(0) = A$. The above formulation \eqref{eq:hhcharac} characterizes, apart from the complex shifted Laplace (CSL) operator, also the concept of complex stretched grid (CSG), where the underlying grid is rotated into the complex plane. This results in a damped Helmholtz problem that is equivalent to complex shifted Laplacian, see \cite{reps2010indefinite}. For the remainder of this text we however assume $P(\beta,\bold{x}) = -\beta i k^2(\bold{x})$ as suggested by \eqref{eq:csl_sys}.

We define an operator functional $f$ based on the general shifted Laplacian operator $M$ as follows:
\begin{equation}
f(\beta) := M(\beta)^{-1} = (-\Delta - (1+\beta i) k^2)^{-1}.
\end{equation}
Choosing $\beta \equiv 0$ in the above expression results in the inverse of the original Helmholtz operator $A = M(0)$, whereas choosing $\beta > 0$ yields the inverse of the shifted Laplace operator $M(\beta)$. The derivatives of the functional $f$ are given by 
\begin{equation} \label{eq:deriv}
f^{(n)}(\beta) = n! \, (k^2 i)^n \, (-\Delta - (1+\beta i) k^2)^{-(n+1)},
\end{equation}
for any $n \in \mathbb{N}$. Constructing a Taylor series expansion \cite{morse1953methods,bhatia1997matrix,arfken2011mathematical} of $f(\beta)$ around a fixed shift $\beta_0 \in \mathbb{R}^{+}$ leads now to the following expression
\begin{equation} \label{eq:series}
f(\beta) = \sum_{n=0}^\infty \frac{f^{(n)}(\beta_0)}{n!}(\beta-\beta_0)^n,
\end{equation}
where the derivatives $f^{(n)}(\beta_0)$ are defined by \eqref{eq:deriv}. Note that the derivatives of $f$ in \eqref{eq:series} are negative powers of the complex shifted Laplace operator $M(\beta_0)$. By evaluating the functional $f(\beta)$ in $\beta = 0$ and by choosing a sufficiently large positive value for $\beta_0$, equation \eqref{eq:series} yields an approximation of the original Helmholtz operator inverse in terms of CSL operator inverses, i.e.
\begin{align} \label{eq:series_orig}
f(0) = M(0)^{-1}	&= \sum_{n=0}^\infty \, (-\beta_0)^n \, \frac{f^{(n)}(\beta_0)}{n!}   \notag \\ 
									&= \sum_{n=0}^\infty \, (-\beta_0 k^2 i)^n \, (-\Delta - (1+\beta_0 i) k^2)^{-(n+1)} \notag \\
									&= \sum_{n=0}^\infty \, (-\beta_0 k^2 i)^n \, M(\beta_0)^{-(n+1)}.
\end{align}
Hence, the computation of the infinite series of easy-to-compute inverse CSL operators $M(\beta_0)$ with an arbitrary shift parameter $\beta_0$ asymptotically results in the exact inversion of the original Helmholtz operator $M(0)= A$. 

\subsubsection{Definition of the expansion preconditioner}

By truncating the expansion in \eqref{eq:series_orig}, we can now define a new class of polynomial Helmholtz preconditioners. For a given $m$, each particular member of this preconditioner class is denoted as the \emph{expansion preconditioner} of degree $m$, 
\begin{equation} \label{eq:EX_def}
EX(m) := \sum_{n=0}^{m-1} \, \alpha_n \, (-\Delta - (1+\beta_0 i) k^2)^{-(n+1)},
\end{equation}
where the coefficients $\alpha_{0},\ldots,\alpha_{m-1}$ are defined as
\begin{equation}
\alpha_n = (-\beta_0 k^2i)^n, \qquad (0 \leq n \leq m-1),
\end{equation} 
by the Taylor series expansion \eqref{eq:deriv}-\eqref{eq:series}. The expansion preconditioner $EX(m)$ is hence a degree $m$ polynomial in the inverse complex shifted Laplace operator $M(\beta_0)^{-1} = (-\Delta - (1+\beta_0 i) k^2)^{-1}$.
The above Taylor series approach appears quite natural. However, other series approximations to the Helmholtz operator inverse may be constructed using alternative choices for the series coefficients. We refer to Section \ref{sec:extensions} for a more elaborate discussion on the choice of the series coefficients. 

\subsubsection{Properties of the expansion preconditioner} \label{sec:properties}

Following the formal definition \eqref{eq:EX_def}, we formulate some essential properties of the $EX(m)$ class of preconditioners in this section. Firstly, one trivially observes that the classic CSL preconditioner is a member of the class of expansion preconditioners. Indeed, the complex shifted Laplace inverse $M(\beta_0)^{-1}$ is the first order term in the Taylor expansion \eqref{eq:series_orig}, and hence we have $M(\beta_0)^{-1} = EX(1)$. 

By including additional terms in the preconditioning polynomial (i.e.~for $m \to \infty$), the $EX(m)$ preconditioner becomes an increasingly accurate approximation to the original Helmholtz operator inverse $A^{-1}$. Hence, the class of $EX(m)$ preconditioners is \emph{asymptotically exact}, since
\begin{equation}\label{eq:limlim}
\lim_{m\to\infty} EX(m) = \sum_{n=0}^\infty \, \alpha_n \, (-\Delta - (1+\beta_0 i) k^2)^{-(n+1)} = M(0)^{-1}.
\end{equation}
This implies that, if we assume that the computational cost of computing the inverse matrix powers in \eqref{eq:limlim} is manageable, $EX(m)$ satisfies \emph{both} conditions (P1) and (P2) for efficient Helmholtz preconditioning suggested in Section \ref{sec:measure}. It should be stressed that (P1) in fact holds asymptotically, and is thus in practice only satisfied when a large number of series terms is taken into account. The $EX(m)$ preconditioner is thus expected to be increasingly more efficient for growing $m$, which suggests a significant reduction in the number of outer Krylov iterations. On the other hand, condition (P2) is satisfied when $m$ is not too large. This creates a trade-off for the value of $m$, which is commented on in Section \ref{sec:costmodel}.

While the approximation precision of the $EX(m)$ preconditioner clearly benefits from the addition of multiple terms in the expansion, note that the accuracy of the $m$-term $EX(m)$ approximation is governed by the truncation error of the series \eqref{eq:series}. This truncation error is of order $\mathcal{O}({\beta_0}^m)$ for any $EX(m)$ preconditioner (with $m>0$), i.e.\
\begin{equation} \label{eq:exporder}
M(0)^{-1} = EX(m) + \mathcal{O}({\beta_0}^m).
\end{equation}
The efficiency of the $EX(m)$ preconditioner is hence also intrinsically dependent on the value of the shift parameter $\beta_0$. However, it is well-known that $\beta_0$ cannot be chosen below a critical value for iterative (multigrid) solver stability, which typically lies around $0.5$ or $0.6$, see \cite{erlangga2006novel,cools2013local}. Consequently, it is clear from \eqref{eq:exporder} that convergence of the Taylor series \eqref{eq:series} is slow, being in the order of ${\beta_0}^m$. This indicates that a large number of terms has to be taken into account in the $EX(m)$ polynomial to obtain a high-precision approximation to the original Helmholtz operator inverse.

\subsubsection{Computational cost of the expansion preconditioner} \label{sec:costmodel}

The inclusion of additional series terms yields a higher-order $EX(m)$ preconditioner polynomial, which is expected to improve preconditioning efficiency as derived above. Therefore, if the computational cost of the CSL inversions would be negligible compared to the cost of applying one Krylov iteration, there would theoretically be no restriction on the number of terms that should be included in $EX(m)$. Unfortunately, even when approximating each CSL inversion by one V-cycle, the computational cost of the CSL inversions is the main bottleneck for the global cost of the solver in practice. Indeed, while the addition of multiple series terms improves performance, it also increases the computational cost of applying the preconditioner. In this section we briefly expound on the computational cost of the $EX(m)$ preconditioner using a simple theoretical cost model.

We model the computational cost of the $EX(m)$ preconditioner by assuming that its cost is directly proportional to the number of CSL operator inversions that need to be performed when solving the preconditioning system. Each additional term in the series \eqref{eq:series} requires exactly one extra shifted Laplace system to be inverted, since the $EX(m)$ polynomial can be constructed as follows:
\begin{align} \label{eq:msystems}
EX(1) w &= \alpha_0 \underbrace{M(\beta_0)^{-1} w}_{:=v_0} , \notag \\ 
EX(2) w &= \alpha_0 v_0 + \alpha_1 \underbrace{M(\beta_0)^{-1} v_0}_{:=v_1} , \notag \\ 
				&\vdots \notag \\
EX(m) w &= \sum_{n=0}^{m-2} \alpha_n v_n + \alpha_{m-1} \underbrace{M(\beta_0)^{-1} v_{m-2}}_{:=v_{m-1}},	
\end{align}
where $w \in \mathbb{R}^N$ is a given vector of size $N$, i.e.~the number of unknowns. Note that all complex shifted Laplace systems in \eqref{eq:msystems} feature the same shift parameter $\beta_0$. Hence, an additional CSL system of the form
\begin{equation} \label{eq:oneofmsystems}
M(\beta_0) \, v_{i} = v_{i-1},  \quad \textrm{with~} v_{-1} := w, \quad (0 \leq i \leq m-1),
\end{equation}
has to be solved for each term in the expansion preconditioner, resulting in a total of $m$ inversions to be performed. We again stress that in practice the shifted Laplace inverse $M(\beta_0)^{-1}$ is never calculated explicitly, but the approximate solution to \eqref{eq:oneofmsystems} is rather computed iteratively by a multigrid V-cycle. 

The question rises whether the reduction in outer Krylov iterations when using the multi-term $EX(m)$ preconditioning polynomial compensates for the rising cost of the additional (approximate) CSL inversions. Let the computational cost of one approximate CSL inversion be denoted as one work unit (1 WU), and let the total computational cost of the $EX(m)$ preconditioner in a complete $EX(m)$-preconditioned Krylov solve be denoted by $\mathcal{C}_{tot}$. If the number of Krylov iterations until convergence (up to a fixed tolerance \texttt{tol}) is $p(m)$, then $\mathcal{C}_{tot} = m \cdot p(m)$ WU. 
Hence, it should hold that $p(m) < C/m$ for some moderate constant $C$ for the cost of the $EX(m)$ preconditioner to support the inclusion of multiple series terms. 
Numerical experiments in Section \ref{sec:experiments} of this work will show that this is generally not the case for the Taylor expansion polynomial in many practical applications, and the classic CSL preconditioner is hence the optimal choice for a preconditioner in the $EX(m)$ class. In the next section we propose several extensions to the $EX(m)$ preconditioner class to further improve its performance.

\section{Extensions and further analysis} 
\label{sec:extensions}

In this section we propose two theoretical extensions to the Taylor series representation \eqref{eq:series} for the inverse Helmholtz operator. These extensions provide essential insights into the expansion preconditioners and aim at further improving preconditioner efficiency. The primary goal is to improve the performance of the expansion preconditioner class, resulting in a more cost-efficient preconditioner with respect to the number of terms $m$. The theoretical results obtained in this section are supported by various numerical experiments in Section \ref{sec:experiments} that substantiate the analysis and illustrate the efficiency of the extended expansion preconditioner.

\subsection{The expansion preconditioner as a stationary iterate}

We first consider an extension of the $EX(m)$ preconditioner class that allows manual optimization of the series coefficients for each degree $m$. To this aim, we illustrate how the $EX(m)$ preconditioner can be interpreted as the $m$-th iterate of a specific fixed-point iteration. Consequently, a class of extended expansion preconditioners is defined by optimizing the fixed-point iteration.

\subsubsection{Taylor series-based polynomial preconditioners as fixed-point iterates}

Recall that the foundation for the Taylor-based $EX(m)$ preconditioner class presented in Section \ref{sec:framework} is the reformulation of the original Helmholtz operator inverse as a Taylor series. Equation \eqref{eq:series_orig} can alternatively be reformulated as
\begin{align} \label{eq:mateq}
f(0) = M(0)^{-1}		&= M(\beta_0)^{-1} \sum_{n=0}^\infty \, (-\beta_0 k^2 i)^n \, M(\beta_0)^{-n}   \notag \\ 
										&= M(\beta_0)^{-1}\, {(I+\beta_0 k^2 i \, M(\beta_0)^{-1})}^{-1},
\end{align}
where the last equation follows from the limit expression 
\begin{equation}
\sum_{n=0}^{\infty} x^n = \frac{1}{1-x} \quad ~ \textrm{for} ~ |x| \leq 1.
\end{equation}
For general (matrix) operators, this series is known in the literature as a Neumann series \cite{werner2006funktionalanalysis}.
Note that for the matrix equation \eqref{eq:mateq} the requirement $|x| \leq 1$ is met if $\rho(-\beta_0 k^2 i \, M(\beta_0)^{-1}) \leq 1$, which is trivially satisfied since $|\lambda_j(M(\beta_0))| \geq |-\beta_0 k^2 i|$ for all $j = 1,\ldots,N$, where we assume $M(\beta_0) \in \mathbb{C}^{N\times N}$. 
For notational convenience, let us denote 
\begin{equation}
L := -\beta_0 k^2 i \, M(\beta_0)^{-1}.
\end{equation}
so that the last line in \eqref{eq:mateq} reads
\begin{equation} \label{eq:mateq2}
M(0)^{-1} = M(\beta_0)^{-1}\, (I-L)^{-1}.
\end{equation}
Equation \eqref{eq:mateq2} shows that inverting the indefinite Helmholtz operator $M(0)$ is equivalent to subsequently inverting the operator $M(\beta_0)$ followed by the inversion of the operator $(I-L)$. The first operator is simply the inverse of a CSL operator and can easily be solved iteratively. However, the second inversion is non-trivial, as it requires the solution $u$ to the system
\begin{equation} \label{eq:lin_sys}
(I-L)u = b,
\end{equation}
given a right-hand side $b \in \mathbb{R}^N$. The linear system \eqref{eq:lin_sys} can alternatively be formulated as a fixed-point iteration (or stationary iterative method)
\begin{equation} \label{eq:iter_jac}
u^{(m+1)} = L u^{(m)} + b, \qquad m > 0.
\end{equation}
One observes that by setting $u^{(0)}=0$, the $m$-th iterate of this fixed-point method generates the $EX(m)$ preconditioner, since
\begin{equation}
u^{(m)} = \left(\sum_{n=0}^{m-1} L^n\right) b \approx (I-L)^{-1} \, b,
\end{equation}
which implies
\begin{equation} \label{eq:pol_taylor}
M(0)^{-1} \, b  \approx M(\beta_0)^{-1} \, u^{(m)} = \left( M(\beta_0)^{-1} \sum_{n=0}^{m-1} L^n \right) b = EX(m) \, b.
\end{equation}
The fixed-point iteration \eqref{eq:iter_jac} thus asymptotically generates the Taylor series expansion \eqref{eq:series}. 

The truncation analysis in Section \ref{sec:properties} indicated that the Taylor series displays a slow convergence. Alternatively, convergence behavior can now be analyzed by studying the convergence of the fixed-point iteration \eqref{eq:iter_jac}, which is governed by the spectral radius of the iteration matrix 
\begin{equation}
\rho(L) = \rho\left(-\beta_0 i k^2 \, (-\Delta - (1+\beta_0 i) k^2)^{-1}\right).
\end{equation}
For $\beta_0 > 0$ this spectral radius tends to be relatively close to one, since
\begin{equation}
\rho(L) = \max_{1 \leq j \leq N} \left| \frac{ -\beta_0 i k^2 }{ \lambda_j - (1+\beta_0 i) k^2 } \right| = \left(\min_{1 \leq j \leq N} \left| 1 -\frac{ \lambda_j-k^2}{ \beta_0 i k^2 } \right|\right)^{-1} \approx 1,
\end{equation}
where $\lambda_j$ $(1\leq j\leq N)$ are the eigenvalues of the negative Laplacian. Hence, the slow convergence of the Taylor series is apparent from the spectral properties of the fixed-point iteration. 

\subsubsection{Weighted fixed-point iteration to improve convergence}

To obtain a convergence speed-up, the fixed-point iteration \eqref{eq:iter_jac} can be substituted by a more general weighted stationary iteration
\begin{equation} 
u^{(m+1)} = (1-\omega) u^{(m)} + \omega L u^{(m)} + \omega b, \qquad \omega \in [0,2], \quad m > 0,
\end{equation}
for $b \in \mathbb{R}^N$, which can alternatively be written as
\begin{equation} \label{eq:iter_wjac}
u^{(m+1)} = \tilde{L} u^{(m)} + \omega b, \qquad \omega \in [0,2], \quad m > 0,
\end{equation}
using the notation $\tilde{L} := (1-\omega) I + \omega L$. Setting $u^{(0)} = 0$ as the initial guess, this iteration constructs a different class of polynomial expansion preconditioners for any choice of $\omega \in [0,2]$, as follows;
\begin{equation}\label{eq:gen_ex_prec}
EX_{\omega}(m) \, b := M(\beta_0)^{-1} \, u^{(m)},  \quad m > 0,
\end{equation}
where $u^{(m)}$ is given by \eqref{eq:iter_wjac}. We call this class of preconditioners the \emph{extended expansion preconditioner} of degree $m$, and denote them by $EX_{\omega}(m)$ to indicate their dependency on $\omega$.

The parameter $\omega$ allows us to modify the coefficients of the series expansion to obtain a more suitable truncated series approximation to the original Helmholtz inverse. Note that the Taylor expansion preconditioner $EX(m)$ can be constructed from iteration \eqref{eq:iter_wjac} by setting the parameter $\omega = 1$. Additionally, note that choosing $\omega = 0$ trivially yields the CSL preconditioner $M(\beta_0)^{-1} = EX_{0}(m)$ for all $m > 0$.

The careful choice of the parameter $\omega \in [0,2]$ in \eqref{eq:gen_ex_prec} possibly results in a series that converges faster than the Taylor series generated by \eqref{eq:iter_jac}. Indeed, the parameter $\omega$ can be chosen to modify the polynomial coefficients such that 
\begin{equation} 
\rho(\tilde{L}) < \rho(L),
\end{equation}
yielding a series that converges faster than the Taylor series \eqref{eq:series}. Consequently, for the right choice of $\omega$, the $EX_{\omega}(m)$ truncated series results in a more efficient preconditioner than the original $EX(m)$ polynomial of the same degree. We refer to the numerical results in Section \ref{sec:experiments} to support this claim.

\subsection{GMRES-based construction of the expansion polynomial} 

The optimization of the series coefficients through the choice of the parameter $\omega$ in the $EX_{\omega}(m)$ operator can be generalized even further by letting the coefficients of the series expansion vary freely. Moreover, the coefficients can be optimized depending on the degree $m$ of the preconditioning polynomial. By replacing the stationary fixed-point iterations in \eqref{eq:iter_jac}-\eqref{eq:iter_wjac} by a more advanced Krylov solution method, an optimal degree $m$ polynomial approximation to the original Helmholtz operator can be constructed.

\subsubsection{Optimization of the expansion preconditioner} \label{sec:gmres_poly}

For the fixed-point method \eqref{eq:iter_jac} we essentially constructed the Taylor polynomial $EX(m)$ using a fixed linear combination of the following basis polynomials 
\begin{equation}
\mathcal{R}(m) = \left\{ L, \, L \,(I + L), \, L \left(I + L + L^2\right), \, \ldots \right\}.
\end{equation}
Alternatively, the extended expansion preconditioner $EX_{\omega}(m)$ was formed as a fixed linear combination of the basis polynomials 
\begin{equation}
\mathcal{S}(m) = \left\{ \omega L, \, (2\omega-\omega^2) L + \omega^2 L^2, \, \ldots \right\},
\end{equation}
for the weighted fixed-point iteration \eqref{eq:iter_wjac}. Note that in this section $L$ designates the inverse of the CSL operator as before, up to scaling by a scalar constant, i.e.
\begin{equation}
L := M(\beta_0)^{-1}.
\end{equation}

As a direct generalization of the above constructions, we now consider the coefficients in each step of the iterative procedure to be variable. This boils down to constructing the preconditioning polynomial from the monomial basis
\begin{equation} \label{eq:mon_basis}
\mathcal{T}(m) = \spann\left\{ L, \, L^2, \, L^3, \, \ldots, \, L^m \right\}.
\end{equation}
Since the preconditioning polynomial asymptotically results in the exact Helmholtz operator inverse, we can alternatively solve the preconditioning system 
\begin{equation}
L \, v = g, \qquad \textrm{with~~} v = A u \textrm{~~and~~} g = L f,
\end{equation} 
using $m$ steps of GMRES \cite{saad1986gmres}, which results in construction of an $m$-term polynomial from the Krylov basis
\begin{equation} \label{eq:kryl_L}
\mathcal{K}_m(L, \, A \, r_0) = \spann\left\{ A \, r_0, \, L \, A \, r_0, \, L^2 \, A \, r_0, \, \ldots, \, L^{m-1} \, A \, r_0 \right\}.
\end{equation}
After an additional multiplication with the operator $L$, the $m$-th Krylov subspace exactly generates a preconditioning polynomial from the basis $\mathcal{T}(m)$ \eqref{eq:mon_basis}. Hence, a generalized expansion preconditioner can be constructed by applying $m$ steps of GMRES on the system $L v = g$, which allows for a free choice of the polynomial coefficients. Moreover, since GMRES minimizes the residual over the $m$-th Krylov subspace, the resulting $m$-term preconditioner is the optimal polynomial approximation of degree $m$ to the exact Helmholtz inverse. These concepts resemble the principles of polynomial smoothing by a GMRES(m)-based construction, see \cite{calandra2012two}.

\subsubsection{Simultaneous construction of preconditioner and Krylov solver basis}

Further extending the above methodology, we outline the theoretical framework for an integrated construction of the preconditioner polynomial in the Krylov subspace construction at the solver runtime level. The key notions in this section show some similarities to the work on multi-preconditioned GMRES in \cite{greif2011multi}. We additionally refer to the closely related literature on flexible Krylov solvers \cite{saad1993fgmres,simoncini2003flexible}. 

Consider the Krylov method solution to the $EX(m)$-preconditioned Helmholtz system 
\begin{equation}
EX(m) \, A u = EX(m) \, f,
\end{equation}
for a fixed polynomial degree $m$. Using GMRES to solve this system, the $k$-th residual $r_k = EX(m) (f-Au^{(k)})$ is minimized over the Krylov subspace
\begin{equation} \label{eq:kryl_EX}
\mathcal{K}_k(EX(m) \, A, r_0) = \spann\left\{ r_0, \, EX(m) \, A \, r_0, \, \ldots, \, (EX(m) \, A)^{k-1} \, r_0 \right\}.
\end{equation}
Note that this basis spans an entirely different subspace compared to the Krylov subspace \eqref{eq:kryl_L}. Indeed, for the basis terms in \eqref{eq:kryl_EX} the preconditioning polynomial degree is fixed at $m$ while the power of the Helmholtz operator $A$ is variable. To form the generalized $EX(m)$ polynomial in \eqref{eq:kryl_L} on the other hand, the powers of the CSL inverse $L$ vary but the power of $A$ is fixed at one. 

A combination of the two principles characterized by \eqref{eq:kryl_L} and \eqref{eq:kryl_EX} can be made by embedding the iterative procedure for the construction of the expansion preconditioner polynomial $EX(m)$ in the governing Krylov solver. This results in a mixed basis consisting of a structured mixture of powers of the inverse CSL operator $L$ and the original system matrix $A$ applied to the initial residual $r_0$. We denote the mixed basis corresponding to the $EX(m)$ preconditioner by
\begin{equation} \label{eq:mixed_base}
\mathcal{K}_k^{EX(m)}(A,r_0) = \spann\left\{ L^{i \cdot j} A^j \, r_0 \, :  1 \leq i < m, ~ 0 \leq j < k \right\},
\end{equation}
for any $m \geq 1$. One trivially observes from the definition \eqref{eq:mixed_base} that
\begin{equation}
\mathcal{K}_k^{EX(m-1)}(A,r_0) \subset \mathcal{K}_k^{EX(m)}(A,r_0),
\end{equation}
which generalizes the embedding of the $EX(1)$ or CSL preconditioner in the class of $EX(m)$ expansion operators to the mixed basis setting. 
The subspace spanned by \eqref{eq:mixed_base} theoretically allows for the simultaneous construction of the preconditioning polynomial and the solution of the preconditioned system. However, the mixed basis $\mathcal{K}_k^{EX(m)}(A,r_0)$ generally does not span a Krylov subspace for any $m > 1$, making its practical construction non-trivial.

As mentioned earlier, the addition of extra terms in the $EX(m)$ polynomial improves the polynomial approximation to the exact Helmholtz operator inverse, resulting in faster convergence in terms of outer Krylov iterations. This implies lower powers of the Helmholtz operator $A$ in the mixed basis \eqref{eq:mixed_base}. However, the addition of extra terms in the polynomial $EX(m)$ also increases the number of vectors constituting the mixed basis, and hence gives rise to a higher computational cost of the total method. This trade-off between preconditioner approximation precision and computational cost in function of the number of terms $m$ is apparent from \eqref{eq:mixed_base}.

As a final remark, note that the extensions proposed in this section are mainly intended as an insightful theoretical framework. In practice the GMRES-based extended preconditioner proposed in Section \ref{sec:gmres_poly} is unlikely to perform significantly better than the extended $EX_{\omega}(m)$ preconditioner proposed in the previous section. This is a consequence of the fact that, given a sufficiently large shift parameter $\beta$, the convergence speed of any monomial-based series of this type is slow, as was already pointed out in Section \ref{sec:properties}.

\section{Numerical results} 
\label{sec:experiments}

In this section we present experimental results that illustrate the practical application of the class of expansion preconditioners to enhance the Krylov convergence on a 1D and 2D Helmholtz benchmark problem. The primary aim is to validate the $EX(m)$ expansion preconditioner for degrees $m > 1$ and illustrate the asymptotic behavior of the $EX(m)$ preconditioner as $m \to \infty$. Initial numerical experiments in Sections \ref{sec:1dhhprob}-\ref{sec:1dperfanal} will use exact inverses of the complex shifted Laplace operators appearing in the polynomial preconditioner. We then introduce a multigrid V(1,1)-cycle as an approximate solver for the CSL systems in the expansion. The performance of the $EX(m)$ preconditioner is consequently compared to that of the classic complex shifted Laplace or $EX(1)$ preconditioner. Additionally, the extensions to the class of generalized $EX_{\omega}(m)$ preconditioners and the combination of the preconditioner polynomial and Krylov basis construction (see Section \ref{sec:extensions}) are shown to display the potential to improve the preconditioner's efficiency.

\subsection{Problem setting: a 1D constant wavenumber Helmholtz problem with absorbing boundary conditions}\label{sec:1dhhprob}

Consider the one-dimensional constant wavenumber Helmholtz model problem on the unit domain
\begin{equation} \label{eq:model_problem}
( -\Delta - k^2 ) \, u(x) = f(x), \quad x \in \Omega = [0,1],
\end{equation}
where the right-hand side $f(x)$ represents a unit source in the domain center. The wavenumber is chosen to be $k^2 = 2 \times 10^4$. The equation is discretized using a standard Shortley-Weller finite difference discretization \cite{shortley1938numerical}, required to treat the absorbing boundary layers (see below). The unit domain $\Omega$ is represented by an $N + 1 = 257$ equidistant point grid, defined as $\Omega^h = \{x_j = j h, 0 \leq j \leq N\}$, respecting the physical wavenumber criterion $kh = 0.5524 < 0.625$ for a minimum of 10 grid points per wavelength \cite{bayliss1985accuracy}. 

To simulate outgoing waves near the edges of the numerical domains, we use \emph{exterior complex scaling} \cite{simon1979definition,chew19943d}, or ECS for short, adding absorbing layers to both sides of the numerical domain. The absorbing layers are implemented by the addition of two artificial complex-valued extensions to the left and right of the domain $\Omega$, defined by the complex grid points $\{z_j = \exp(i \theta_{ECS}) \, x_j\}$, where $x_j = j h$, for $-N/4 \leq j < 0$ and $N < j \leq 5/4 N$. The ECS complex scaling angle that determines the inclination of the extensions in the complex plane is chosen as $\theta_{ECS} = \pi/6$. The two complex-valued extensions feature $N/4$ grid points each, implying the discretized Helmholtz equation takes the form of an extended linear system
\begin{equation}\label{eq:num_exp_hh_sys}
Au = f,
\end{equation}
where $A \in \mathbb{C}^{\frac{3}{2} N \times \frac{3}{2} N}$. 
The discretized right-hand side $f = (f_j) \in \mathbb{C}^{\frac{3}{2}N}$ is defined as 
\begin{equation}
f_j = f(x_j) = 	\left\{
		\begin{array}{ll}
		1  & \mbox{for }  j = N/2, \\
		0  & \mbox{elsewhere,}
		\end{array} 
\right.
\end{equation}
representing a unit source located in the center of the domain for this example.

\begin{figure}[t!] 
\begin{center}
\includegraphics[width=0.48\textwidth]{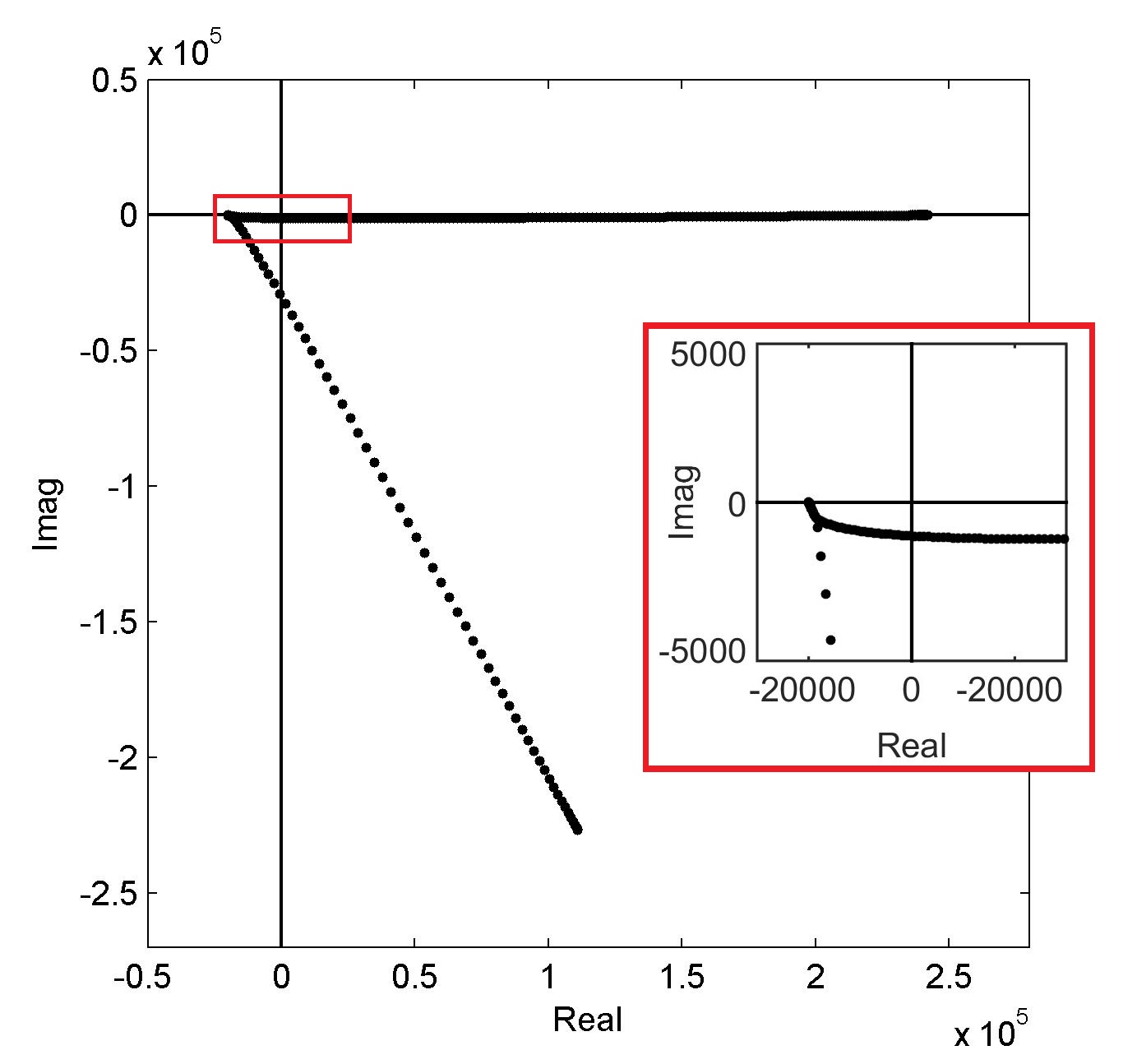}
\includegraphics[width=0.48\textwidth]{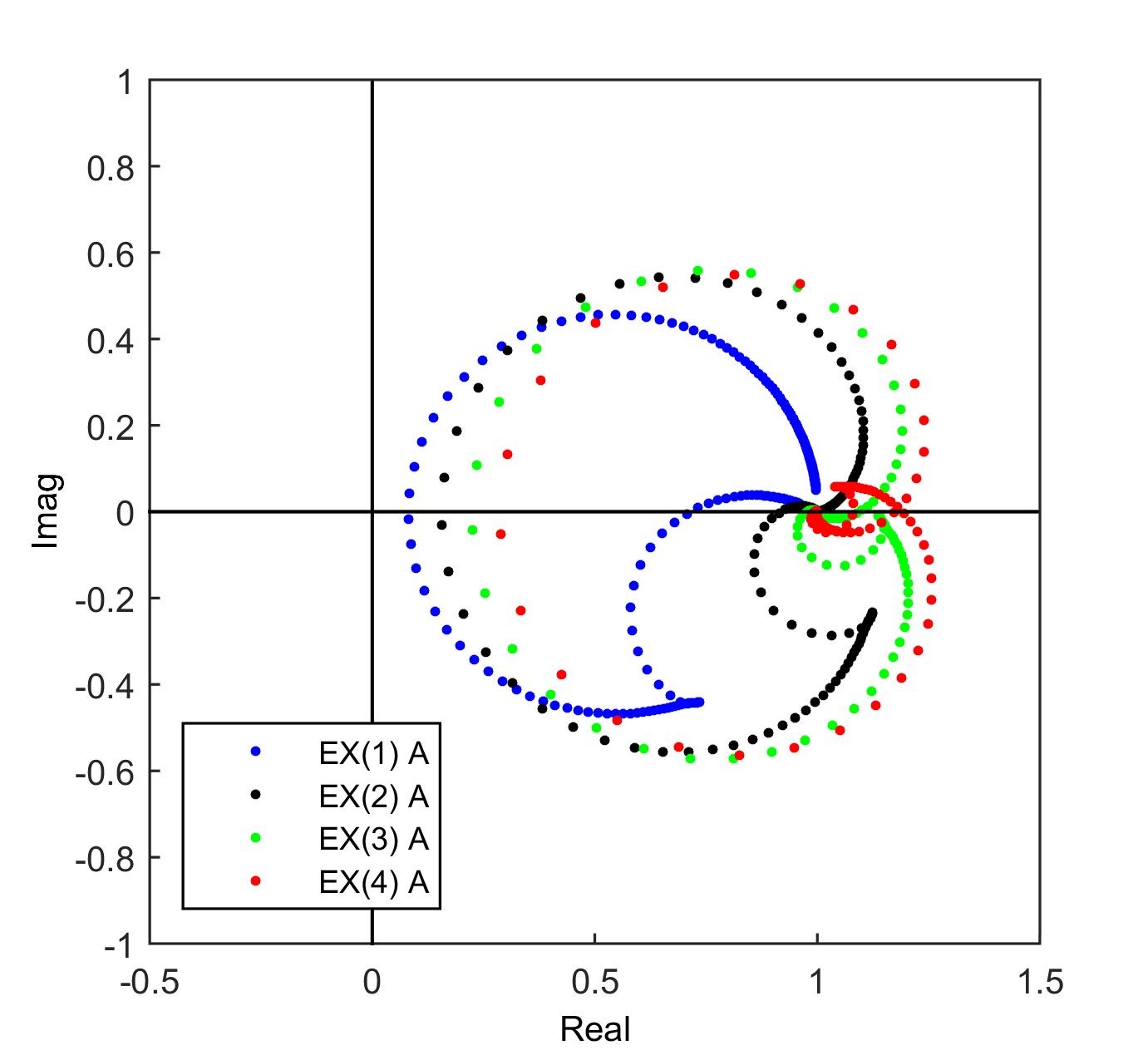}
\caption{Spectral analysis of the discretized 1D Helmholtz model problem \eqref{eq:model_problem}. Exact preconditioner inversion. Left: spectrum of the Helmholtz operator $A$ with ECS absorbing boundary conditions. Right: spectrum of the polynomial preconditioned operator $EX(m)A$ for various values of $m$. The spectrum becomes more clustered around $1$ for increasing values of $m$.}
\label{fig:EXspectra}
\end{center}
\end{figure}

The Helmholtz model problem \eqref{eq:model_problem} is solved using $EX(m)$-pre\-con\-di\-ti\-o\-ned BiCGStab \cite{van1992bicg} up to a relative residual tolerance $\|r_p\|/\|r_0\| < \textrm{\texttt{tol}} = 1$e$-8$.  The CSL operator inverses in the $EX(m)$ polynomial are either computed exactly using LU factorization for the purpose of analysis, or approximated using one multigrid V-cycle as is common in realistic applications. Note that the complex shift parameter $\beta$ in the CSL operators is chosen as $\beta = 0.6$, which guarantees multigrid V(1,1)-cycle stability \cite{cools2013local}.

\subsection{Spectral analysis of the expansion preconditioner}\label{sec:1dspecanal}

To analyze the efficiency of the $EX(m)$ Helmholtz preconditioner we perform a classic eigenvalue analysis of the $EX(m)$-preconditioned Helmholtz operator. For convenience of analysis, the CSL inversions in the $EX(m)$ preconditioning polynomial are solved using a direct method in this section.

The typical pitchfork shaped spectrum of the indefinite Helmholtz operator with ECS boundary conditions is shown in the left panel of Figure \ref{fig:EXspectra}. The leftmost eigenvalue is located near $-k^2 = -2 \times 10^4$, while the rightmost eigenvalue is close to $4/h^2 - k^2 \approx 2.4 \times 10^5$, cf.\ \cite{reps2012analyzing}. The right panel of Figure \ref{fig:EXspectra} shows the spectrum of the $EX(m)$-preconditioned Helmholtz operator for various numbers of terms in the Taylor polynomial $EX(m)$. Note how the spectra become more clustered around $1$ when additional series terms are taken into account, illustrating the asymptotic exactness of the $EX(m)$ preconditioner class. 

\begin{figure}[t!] 
\begin{center}
\includegraphics[width=0.48\textwidth]{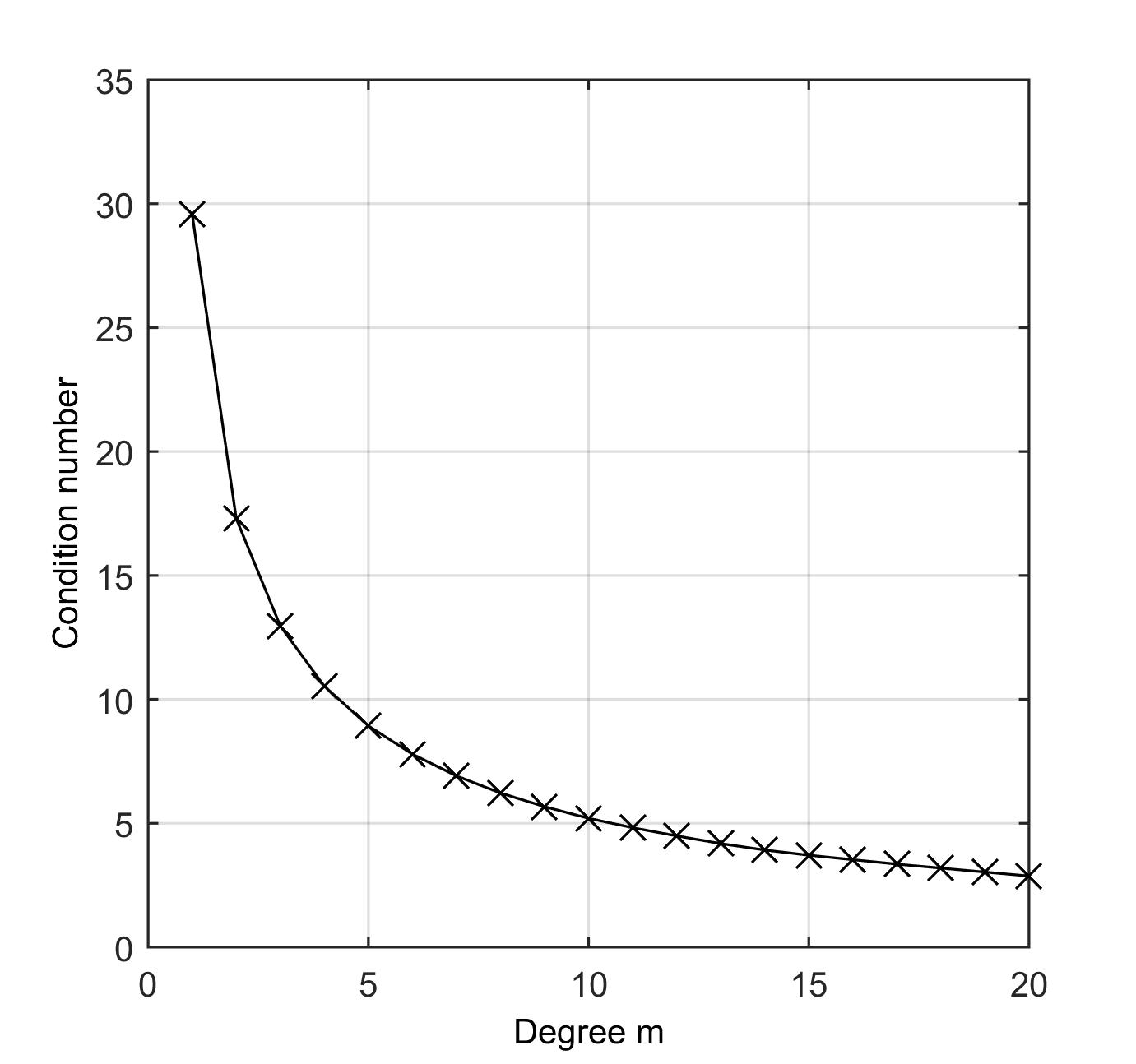}
\includegraphics[width=0.48\textwidth]{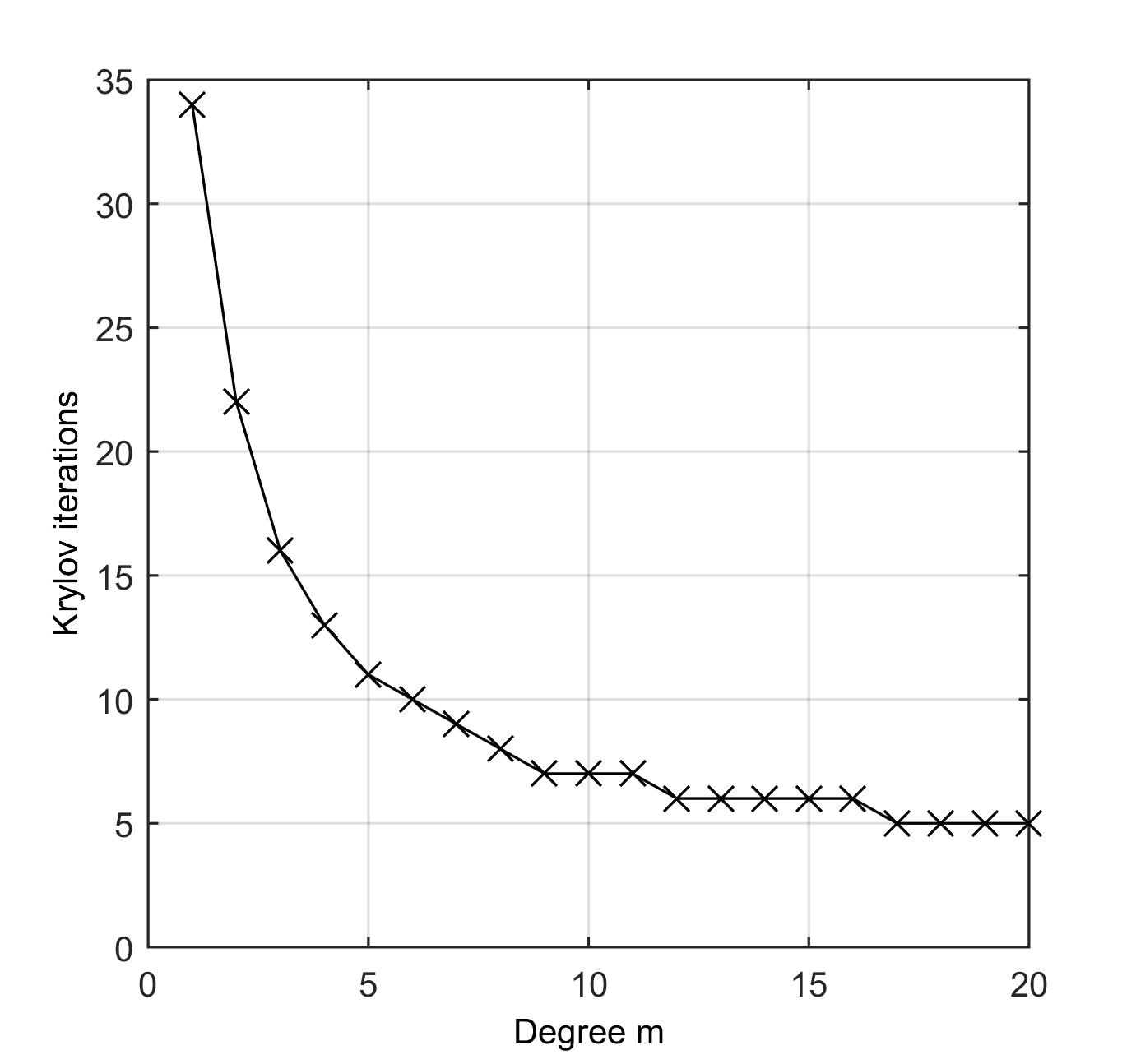}
\caption{Conditioning and $EX(m)$-BiCGStab performance on the discretized 1D Helmholtz model problem \eqref{eq:model_problem}. Exact preconditioner inversion. Left: condition number of the preconditioned operator $EX(m)A$ as a function of $m$. Right: number of $EX(m)$-BiCGStab iterations required to solve the Helmholtz system \eqref{eq:num_exp_hh_sys} as a function of $m$.}
\label{fig:EXcond}
\end{center}
\end{figure}

The condition number of the preconditioned operator $\kappa(EX(m)A)$ is displayed in the left panel of Figure \ref{fig:EXcond} as a function of $m$. One observes that conditioning improves significantly by the addition of extra terms in the polynomial $EX(m)$. This observation is reflected in the number of Krylov iterations required to solve the problem, which is displayed in Figure \ref{fig:EXcond} (right panel) for a range of values of $m$. The number of Krylov iterations (right panel) appears to be directly proportional to the condition number of the preconditioned system (left panel).

\subsection{Performance analysis of the expansion preconditioner}\label{sec:1dperfanal}

It is clear from the spectral analysis that the addition of multiple series terms improves the conditioning of the preconditioned system. In this section, the performance of the $EX(m)$ preconditioner is analyzed using the simple theoretical cost model introduced in Section \ref{sec:costmodel}. 

The experimentally measured performance of the $EX(m)$-BiCGStab solver on the model problem \eqref{eq:model_problem} is displayed in Table \ref{tab:EXperformance2}. The table shows the number of Krylov iterations required to solve system \eqref{eq:num_exp_hh_sys} until convergence up to \texttt{tol} = $1$e$-8$ (first column), the iteration ratio compared to the standard CSL preconditioner (second column), and the effective number of work units (CSL inversions) for the entire run of the method (third column). The $EX(m)$ preconditioner becomes increasingly more efficient in reducing the number of Krylov iterations in function of larger $m$. Comparing e.g.\ the $EX(3)$ preconditioner to the classic $EX(1)$ (CSL) scheme, one observes that the number of Krylov iterations is slightly more than halved. The largest improvement is obtained by adding the first few terms, which is a consequence of the slow Taylor convergence. Note that the computational cost of the Krylov solver itself is not incorporated into this cost model.

\begin{table}[t]
\centering
\begin{tabular}{| C{0.8cm} || C{1.6cm} | C{1.6cm} | C{2.2cm} |}
\hline
$m$ & $p(m)$ & $\frac{p(m)}{p(1)}$ & $m \cdot p(m)$\\
\hline
	1 & 34 & 1.00 & 34 WU \\
	2 & 22 & 0.65 & 44 WU \\
	3 & 16 & 0.47 & 48 WU \\
	4 & 13 & 0.38 & 52 WU \\
	5 & 11 & 0.32 & 55 WU \\
\hline
\end{tabular}
\vspace{0.3cm}
\caption{Performance of $EX(m)$-BiCGStab for different values of $m$ on the discretized 1D Helmholtz model problem \eqref{eq:model_problem}. Exact preconditioner inversion.  Column 1: number of BiCGStab iterations $p(m)$ required to solve the system \eqref{eq:num_exp_hh_sys}. Col.~2: iteration ratio compared to classic CSL. Col.~3: preconditioner computational cost based on $p(m)$.}
\label{tab:EXperformance2}
\end{table}

Although higher-order series approximations clearly result in a qualitatively better preconditioner, the number of Krylov iterations is not reduced sufficiently to compensate for the cost of the extra inversions. Indeed, while the addition of multiple series terms in the $EX(m)$ preconditioner improves the spectral properties of the preconditioned system, the increased computational cost of the extra CSL inversions appears to be a bottleneck for performance. Hence, one observes that standard CSL preconditioning -- which takes only one series term into account -- is the most cost-efficient, requiring a minimum of 34 WU for the entire solve.

\subsection{Multigrid inversion of the expansion preconditioner} \label{sec:multires}

\begin{table}[t]
\centering
\begin{tabular}{| C{0.8cm} || C{1.6cm} | C{1.6cm}  | C{2.2cm} |}
\hline
$m$ & $p(m)$ & $\frac{p(m)}{p(1)}$ & CPU time\\
\hline
1 & 49 & 1.00 & 0.57 s. \\
2 & 39 & 0.80 & 0.68 s. \\
3 & 34 & 0.69 & 0.80 s. \\
4 & 31 & 0.63 & 0.88 s. \\
5 & 30 & 0.61 & 1.02 s. \\
\hline
\end{tabular}
\vspace{0.3cm}
\caption{Performance of $EX(m)$-BiCGStab for different values of $m$ on the discretized 1D Helmholtz model problem \eqref{eq:model_problem}. V(1,1)-cycle approximate preconditioner inversion. Column 1: number of BiCGStab iterations $p(m)$ required to solve the system \eqref{eq:num_exp_hh_sys}. Col.~2: iteration ratio compared to CSL.  Col.~3: CPU time until convergence (in seconds).}
\label{tab:EXmultigrid}
\end{table}

For convenience of analysis a direct inversion of the preconditioning scheme was used in the previous sections. However, in realistic large-scale applications the terms of the $EX(m)$ preconditioner often cannot be computed directly. Instead, the CSL systems comprising the $EX(m)$ polynomial are approximately solved using some iterative method. In this section we use one geometric multigrid V(1,1)-cycle to approximately solve the CSL systems, which is a standard approach in the Helmholtz literature \cite{erlangga2006novel,erlangga2008multilevel}. The V(1,1)-cycle features the traditional linear interpolation and full weighting restriction as intergrid operators, and applies one weighted Jacobi iteration (with parameter 2/3) as a pre- and post-smoother. The choice of the damping parameter $\beta = 0.6$ guarantees stability of the multigrid solver for the inversion of the CSL operators, see \cite{cools2013local}. 

\begin{figure}[t!] 
\begin{center}
\includegraphics[width=0.48\textwidth]{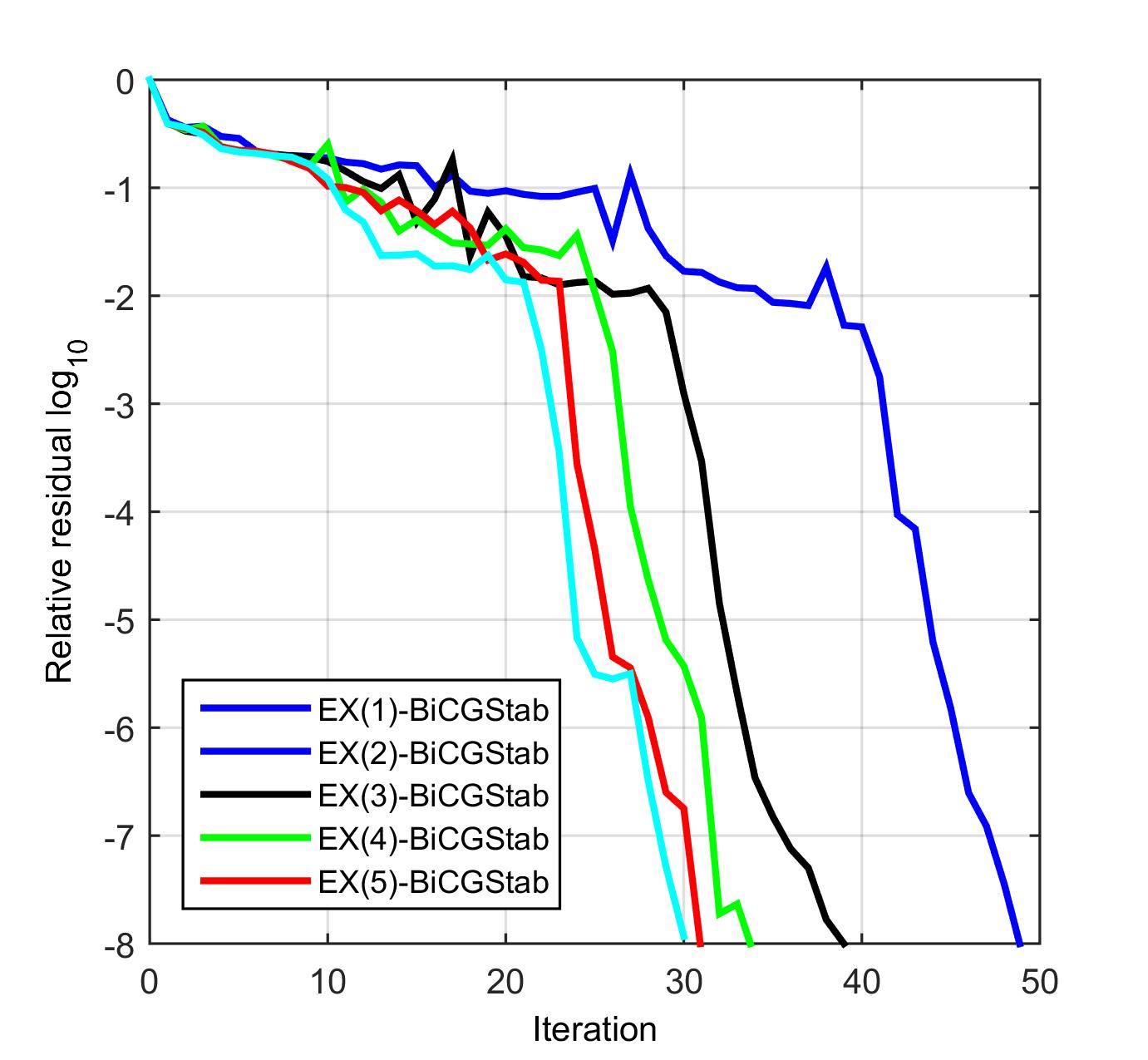}
\includegraphics[width=0.48\textwidth]{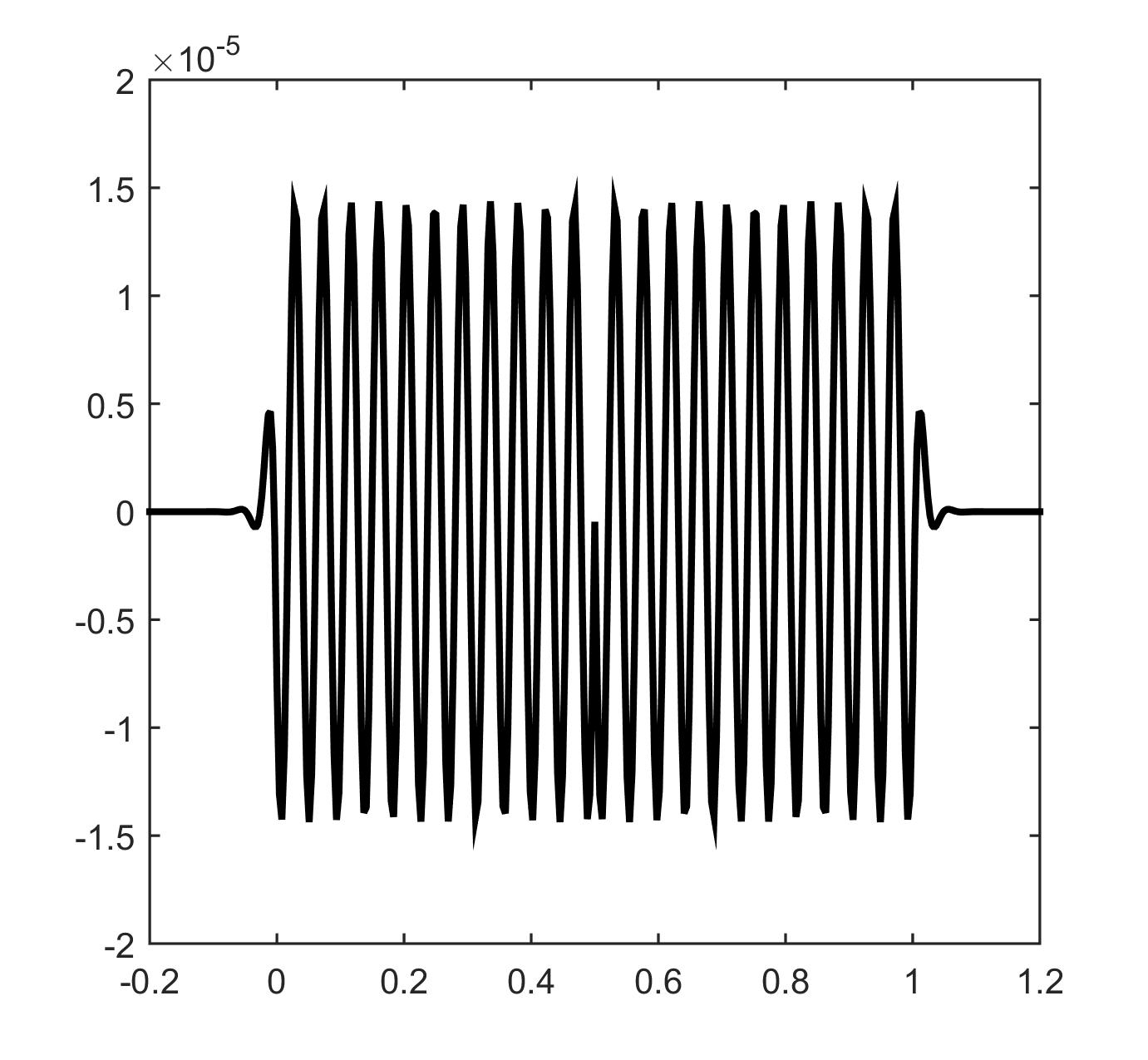}
\caption{Convergence of $EX(m)$-BiCGStab and solution of the discretized 1D Helmholtz model problem \eqref{eq:model_problem}. V(1,1)-cycle approximate preconditioner inversion. Left: $EX(m)$-BiCGStab relative residual history $\|r_p\|/\|r_0\|$ for various values of $m$. Vertical axis in log scale. Right: numerical solution $u(x)$ to the equation \eqref{eq:model_problem} up to the relative residual tolerance $\texttt{tol} = 1$e$-8$. The ECS absorbing boundary layer rapidly damps the solution outside the unit domain $\Omega = [0,1]$.}
\label{fig:EXvcycle}
\end{center}
\end{figure}

Table \ref{tab:EXmultigrid} summarizes the number of Krylov iterations and CPU time\footnote{Hardware specifications: Intel Core i7-2720QM 2.20GHz CPU, 6MB Cache, 8GB RAM. Software specifications: Windows 7 64-bit OS, experiments implemented in MATLAB R2015a.} required to solve system \eqref{eq:num_exp_hh_sys} using $EX(m)$-BiCGStab for different values of $m$. The application of the $EX(m)$ operator is approximately computed using a total of $m$ V(1,1)-cycles, see Section \ref{sec:costmodel}. The corresponding convergence histories are shown in Figure \ref{fig:EXvcycle}. The addition of terms in the $EX(m)$ preconditioner reduces the number of Krylov iterations as expected, although the improvement is less pronounced compared to the results in Table \ref{tab:EXperformance2} due to non-exact inversion of the CSL operators. However, the increased cost to (approximately) compute the additional series terms for larger values of $m$ is clearly reflected in the timings. Hence, in terms of preconditioner computational cost, the classic CSL or $EX(1)$ preconditioner is the most cost-efficient member of the $EX(m)$ preconditioner class for this benchmark problem.

\subsection{Validation of the extended expansion preconditioner}

In this section we validate the generalizations to the expansion preconditioner proposed in Section \ref{sec:extensions} on the 1D Helmholtz model problem \eqref{eq:model_problem}. 

\begin{figure}[t!] 
\begin{center}
\includegraphics[width=0.55\textwidth]{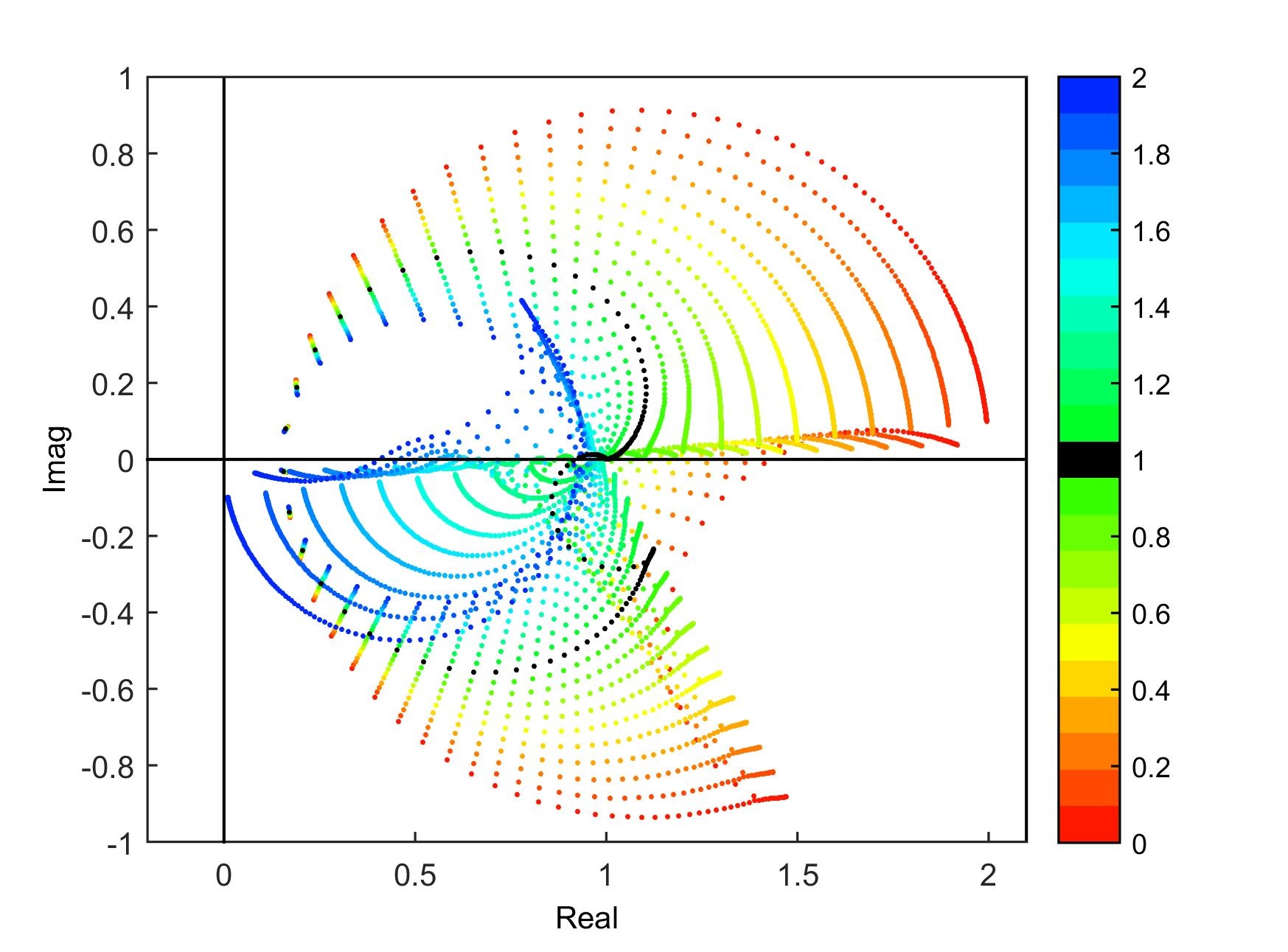}
\includegraphics[width=0.44\textwidth]{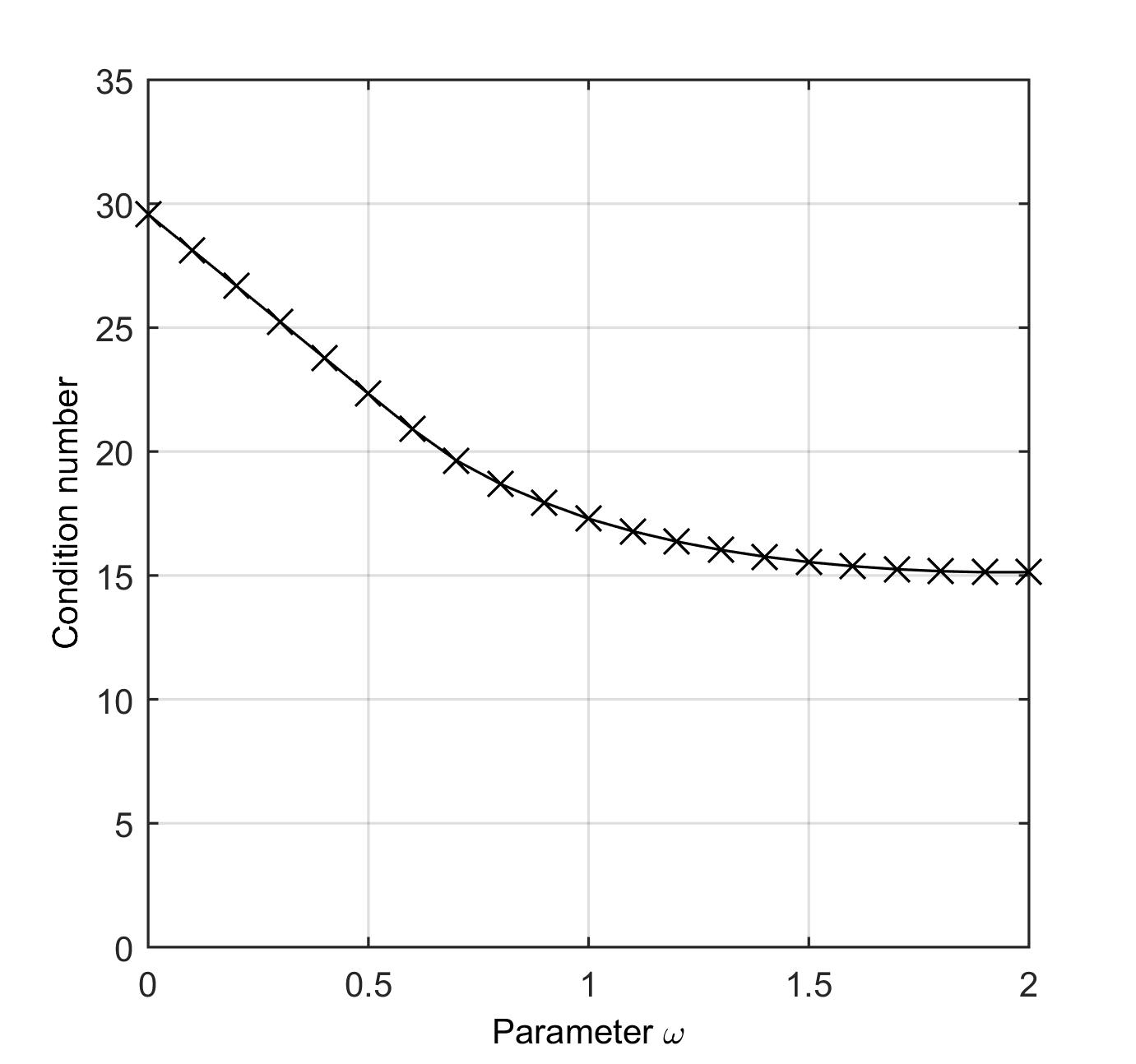}
\caption{Spectral analysis of the discretized 1D Helmholtz model problem \eqref{eq:model_problem}. Exact preconditioner inversion. Left: spectrum of the preconditioned operator $EX_{\omega}(2)A$ for different values of the parameter $\omega$. The spectrum for $EX(m)A$ with $\omega = 1$ is indicated in black, see also Fig.~\ref{fig:EXspectra}. Right: condition number of the preconditioned operator $EX_{\omega}(2)A$ as a function of the weight $\omega$.}
\label{fig:EXomega}
\end{center}
\end{figure}

The weighted fixed-point iteration \eqref{eq:iter_wjac} generates the generalized class of expansion preconditioners $EX_{\omega}(m)$. Figure \ref{fig:EXomega} shows the spectrum (left panel) and condition number (right panel) of the Helmholtz operator preconditioned by the two-term $EX_{\omega}(2)$ polynomial for different values of the parameter $\omega$. Note that the condition number of the standard $EX(2)$ preconditioner ($\omega = 1$) is $\kappa(EX(2)A) = 17.29$. A small improvement in conditioning is achievable through the right choice of the parameter $\omega$, reducing the condition number to $\kappa(EX_{\omega}(2)A) = 15.13$ for parameter choices around $\omega = 2$. With the optimal choice for $\omega$, the condition number of the classic CSL preconditioner $EX_0(2)$ is halved when using $EX_{\omega}(2)$, suggesting a halving of the number of Krylov iterations may be achievable by using $EX_{\omega}(2)$ instead of the classic CSL preconditioner. The smaller condition number implies an increase in performance compared to the $EX(2)$ preconditioner, making the addition of extra terms in $EX_{\omega}(m)$ theoretically cost-efficient.

\begin{figure}[t!] 
\begin{center}
\includegraphics[width=0.55\textwidth]{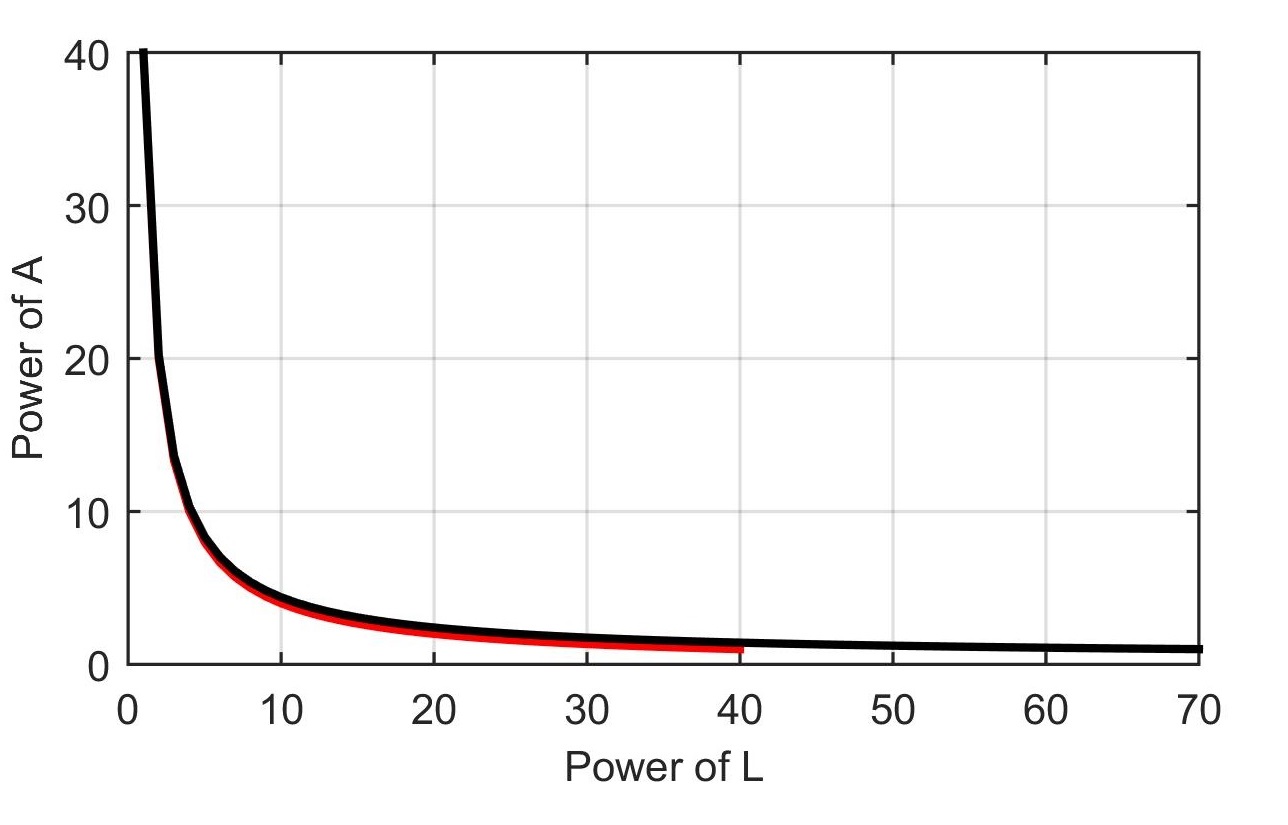}
\caption{Performance of the $EX(m)$ preconditioner formed by the construction of the mixed basis \eqref{eq:mixed_base} on the discretized 1D Helmholtz model problem \eqref{eq:model_problem}.  Maximum power of $L$, i.e.~the polynomial preconditioner degree $m$, vs.\ maximum power of $A$, i.e.~the number of Krylov iterations $p(m)$. Red curve: theoretical upper bound required for cost-efficiency. Black curve: experimentally measured results.}
\label{fig:EXgmres}
\end{center}
\end{figure}

A first step towards a simultaneous construction of the $EX(m)$ preconditioning polynomial and the outer Krylov basis resulting in the mixed basis \eqref{eq:mixed_base} is illustrated in Figure \ref{fig:EXgmres}. The black curve shows the experimentally determined maximum power of $A$ (Krylov iterations) versus the maximum power of $L$ (terms in the preconditioning polynomial) required to solve the Helmholtz benchmark problem \eqref{eq:model_problem} up to a relative residual tolerance \texttt{tol} $= 1$e$-8$. Subject to this tolerance, a solution is found either after 40 $EX(1)$-BiCGStab iterations or alternatively after one $EX(70)$-BiCGStab iteration. Indeed, the incorporation of 70 terms in the $EX(m)$ expansion preconditioner effectively reduces the number of outer Krylov iterations to $1$. However, note that to be cost-efficient with respect to the number of CSL inversions, the same solution should be found using the $EX(40)$ (degree 40) polynomial preconditioner. The red curve represents a constant number of CSL inversions for the total run of the method. To ensure cost-efficiency of the class of $EX(m)$ preconditioners for $m > 1$, the experimental black curve should fall below the theoretical red curve, which is not the case. Hence, the most simple case of the $EX(1)$ or CSL preconditioner can again be considered optimal w.r.t.~cost-efficiency.

\subsection{Problem setting: a 2D constant wavenumber Helmholtz problem with absorbing boundary conditions} 

\begin{figure}[t!] 
\begin{center}
\includegraphics[width=0.48\textwidth]{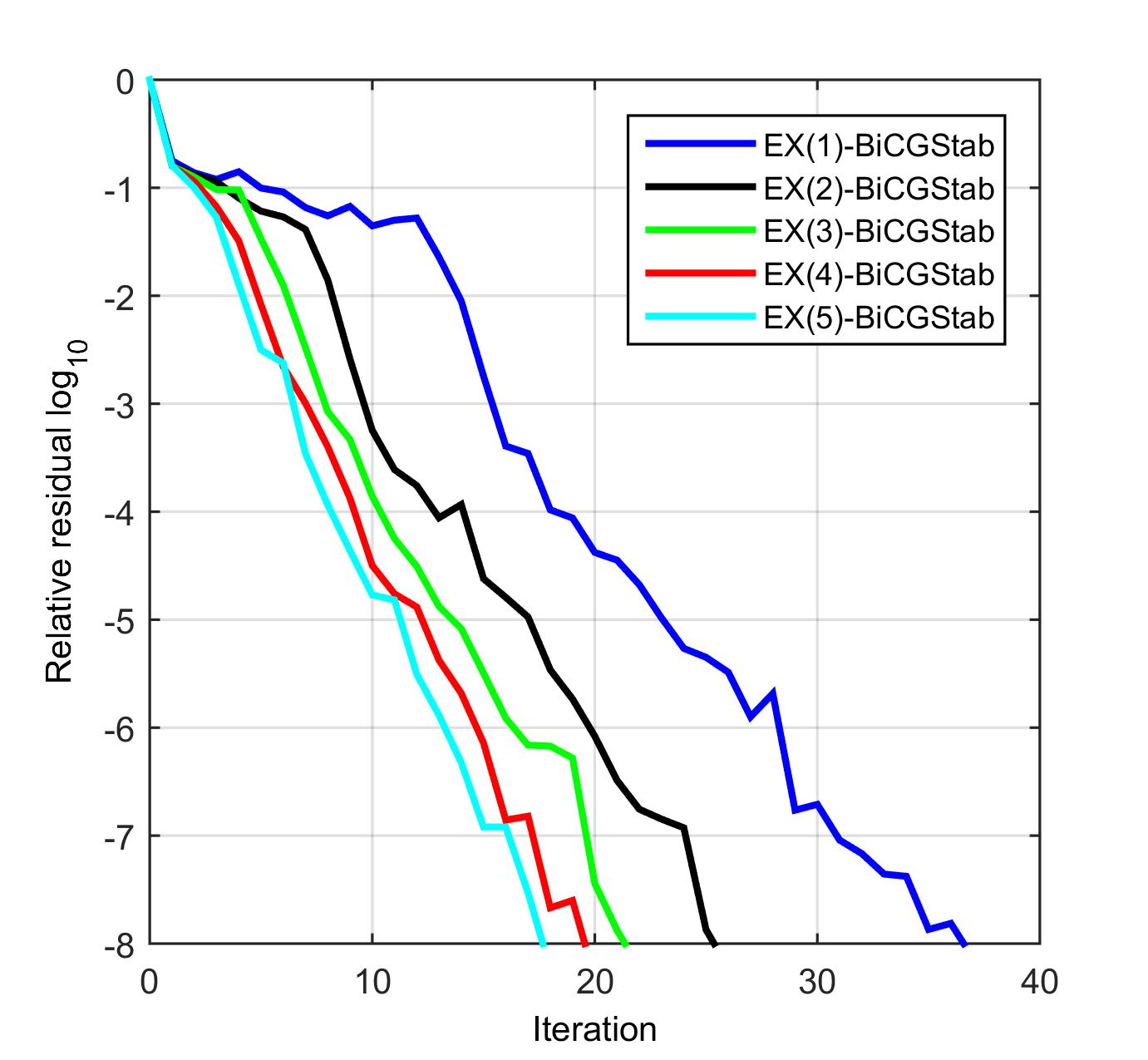}
\includegraphics[width=0.48\textwidth]{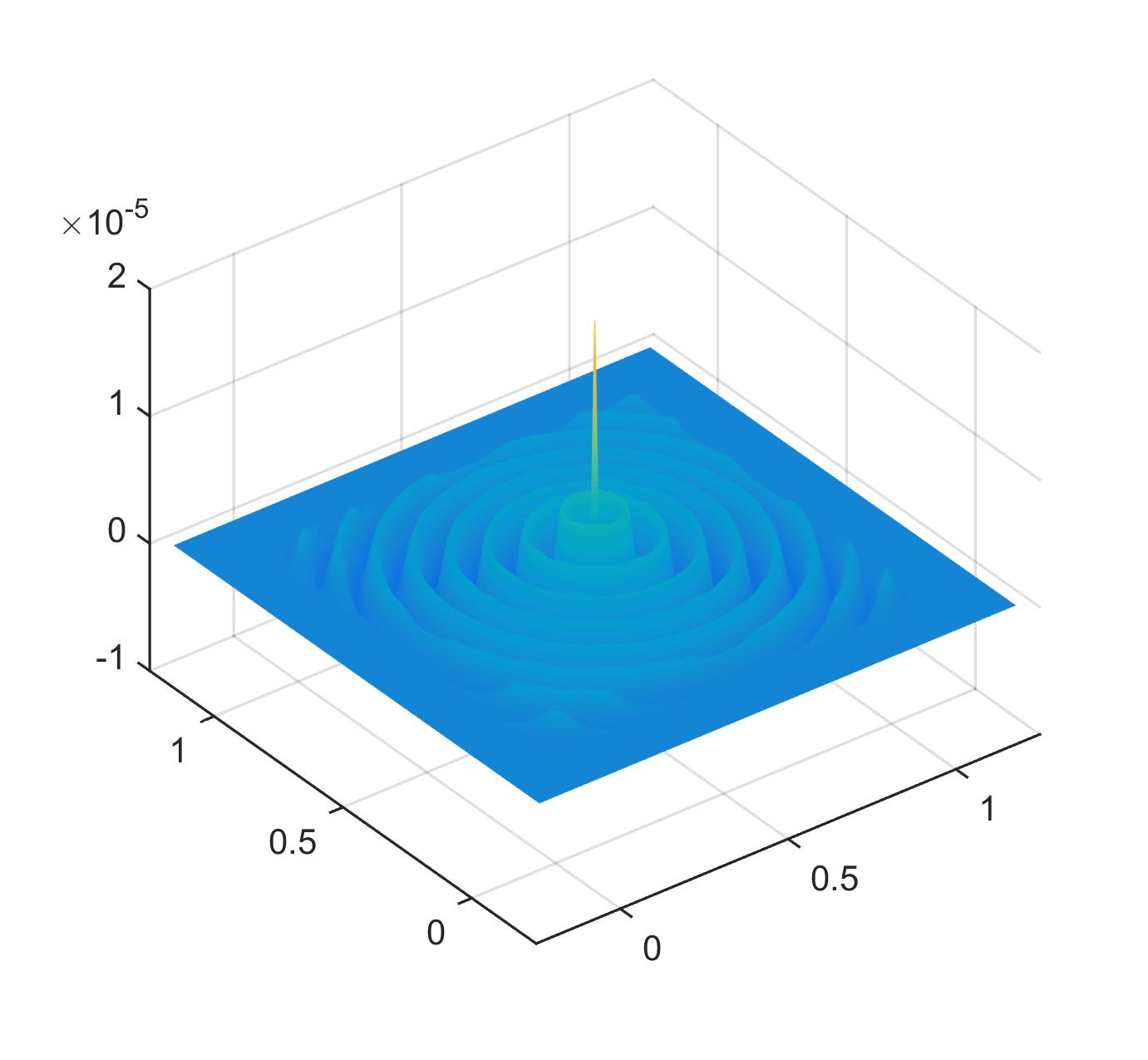} \\
\includegraphics[width=0.48\textwidth]{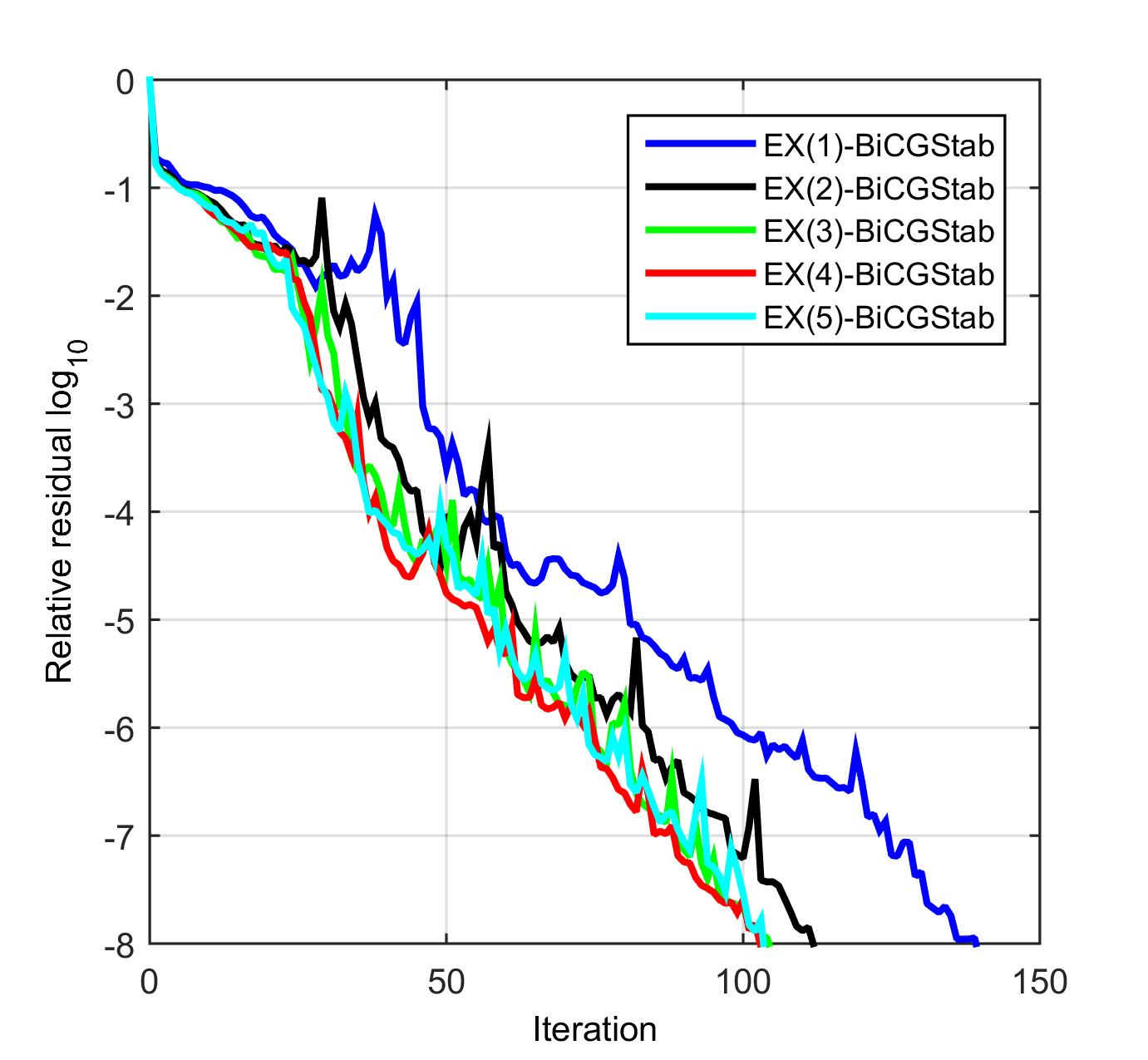}
\includegraphics[width=0.48\textwidth]{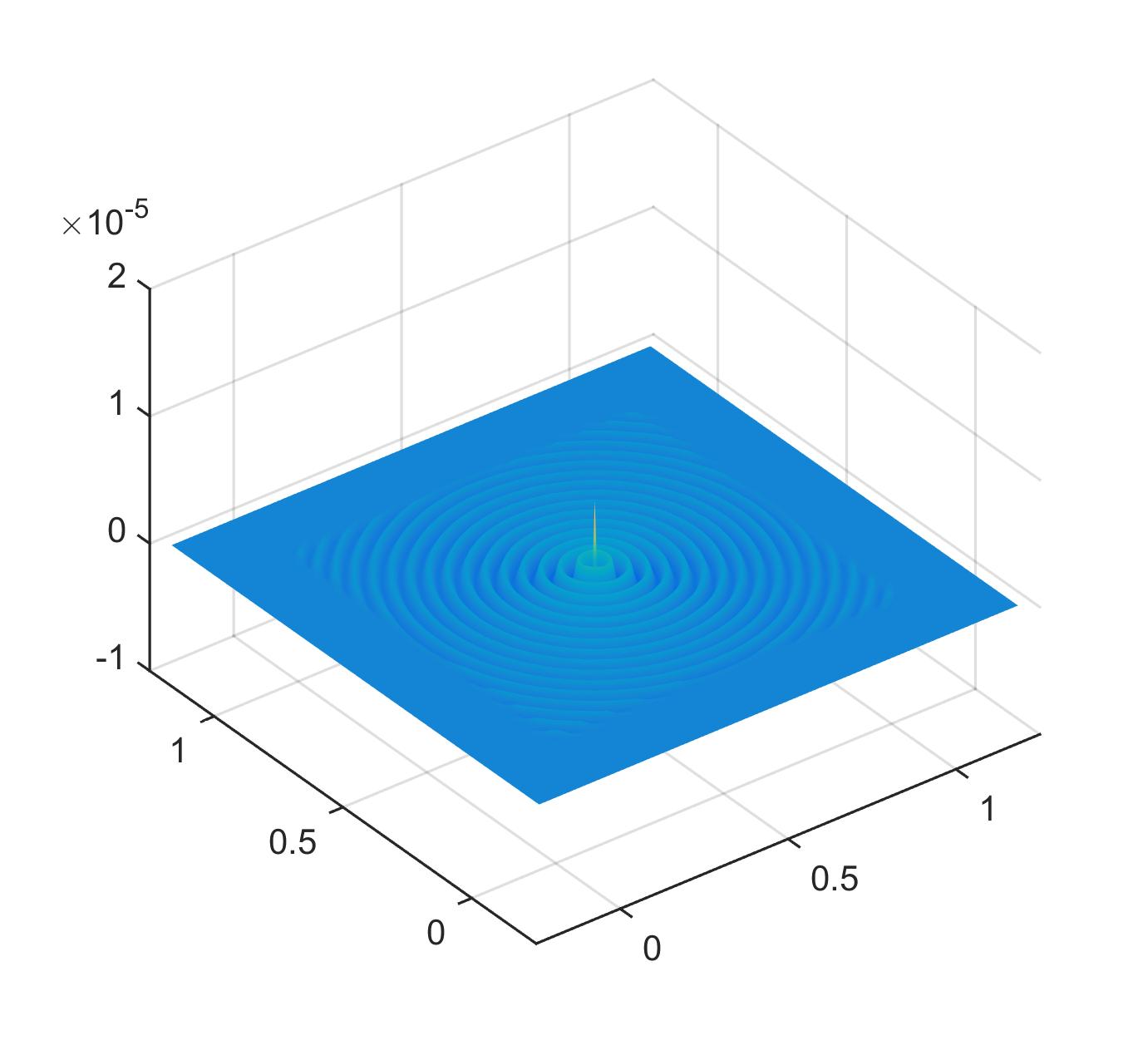}
\caption{Convergence of $EX(m)$-BiCGStab and solutions of the discretized 2D Helmholtz model problem \eqref{eq:model_problem_2d}. Top: wavenumber $k^2 = 5 \times 10^3$ and $N = n_x \times n_y = 128 \times 128$ unknowns. Bottom: wavenumber $k^2 = 2 \times 10^4$ and $N = 256 \times 256$ unknowns. V(1,1)-cycle approximate preconditioner inversion. Left: $EX(m)$-BiCGStab relative residual history $\|r_p\|/\|r_0\|$ for various values of $m$. Vertical axis in log scale. Right: numerical solution $u(x,y)$ to the equation \eqref{eq:model_problem_2d} up to the relative residual tolerance  $\texttt{tol} = 1$e$-8$.}
\label{fig:EXvcycle_2d}
\end{center}
\end{figure}

To conclude this work we extend the above 1D model problem \eqref{eq:model_problem} to a two dimensional Helmholtz problem. Numerical results for solution using $EX(m)$-preconditioned BiCGStab and GMRES are provided, and we comment on the scalability of the expansion preconditioner functionality to higher spatial dimensions. 

Consider the two-dimensional constant wavenumber Helmholtz model problem
\begin{equation} \label{eq:model_problem_2d}
( -\Delta - k^2 ) \, u(x,y) = f(x,y), \quad (x,y) \in \Omega = [0,1]^2,
\end{equation}
where the right-hand side $f(x,y)$ again represents a unit source in the domain center and outgoing wave boundary conditions are implemented using Exterior Complex Scaling with $\theta_{ECS} = \pi/6$. We consider two different wavenumbers, namely $k^2 = 5$e$+3$ and $k^2 = 2$e$+4$, corresponding to a moderate- and high-energetic wave respectively. Equation \eqref{eq:model_problem_2d} is discretized using $n_x = n_y = 128$ (for $k^2 = 5$e$+3$) and $n_x = n_y = 256$ (for $k^2 = 2$e$+4$) real-valued grid points in each spatial dimension, respecting the wavenumber criterion $kh < 0.625$ \cite{bayliss1985accuracy} in every direction. Note that the discretized 2D Helmholtz operator with ECS boundary conditions can be constructed from the 1D Helmholtz operator using Kronecker products, i.e.
$A^{2D} = A_x^{1D} \otimes I_y + I_x \otimes A_y^{1D}$,
where $I_x \in \mathbb{C}^{n_x \times n_x}$ and $I_y \in \mathbb{C}^{n_y \times n_y}$ are identity matrices.

Figure \ref{fig:EXvcycle_2d} shows the $EX(m)$-BiCGStab solver convergence history for various values of $m$ (left) and the solution $u(x,y)$ (right) to the 2D model problem \eqref{eq:model_problem_2d} for different wavenumbers and corresponding discretizations. The $EX(m)$ preconditioner is approximately inverted using $m$ multigrid V(1,1)-cycles with a weighted Jacobi smoother (weighting parameter 4/5). The corresponding number of BiCGStab iterations until convergence up to the relative residual tolerance \texttt{tol}$ = 1$e$-8$ are displayed in Table \ref{tab:EXmultigrid_2d}. The observations from the 1D spectral analysis extend directly to the 2D setting, as the table shows that the use of the $EX(m)$ preconditioner results in a significant reduction of the number of outer Krylov iterations for growing values of $m$.

\begin{table}[t]
\centering {\small
\begin{tabular}{ C{0.8cm} | C{1.4cm} | C{2.0cm} | C{1.4cm} | C{2.0cm} |}
\cline{2-5}
& \multicolumn{2}{c| }{$n_x \times n_y = 128 \times 128$} & \multicolumn{2}{ c| }{$n_x \times n_y = 256 \times 256$} \\
& \multicolumn{2}{c| }{$k^2 = 5$e$+3$} & \multicolumn{2}{ c| }{$k^2 = 2$e$+4$} \\
\hline
\multicolumn{1}{|c||}{$m$} & $p(m)$ & CPU time & $p(m)$ & CPU time\\
\hline
\multicolumn{1}{|c||}{1} & 37 & 14.7 s. & 140 & 157.0 s. \\
\multicolumn{1}{|c||}{2} & 26 & 19.0 s. & 112 & 210.8 s. \\
\multicolumn{1}{|c||}{3} & 22 & 23.1 s. & 105 & 277.9 s. \\
\multicolumn{1}{|c||}{4} & 20 & 26.3 s. & 104 & 351.0 s. \\
\multicolumn{1}{|c||}{5} & 18 & 30.5 s. & 103 & 436.2 s. \\
\hline
& \multicolumn{4}{c| }{EX(m)-preconditioned BiCGStab} \\ 
\cline{2-5}
\end{tabular}
}
\vspace{0.3cm}
\caption{Performance of $EX(m)$-BiCGStab for different values of $m$ on the discretized 2D Helmholtz model problem \eqref{eq:model_problem_2d}. V(1,1)-cycle approximate preconditioner inversion. Column 1 \& 3: number of $EX(m)$-BiCGStab iterations $p(m)$. Col.~2 \& 4: total CPU time until convergence (in seconds).}
\label{tab:EXmultigrid_2d}
\end{table}

Table \ref{tab:EXmultigrid_2d} additionally features the CPU timings for the wavenumbers $k^2 = 5$e$+3$ and $2$e$+4$. Note that although the number of Krylov iterations is reduced as a function of the preconditioner degree $m$, the CPU timings are rising in function of $m$. The computational cost of performing $m$ multigrid V-cycles (compared to just one V-cycle for the CSL preconditioner) has a clear impact on the CPU timings. As a result, it is often advisable in view of cost-efficiency to restrict the expansion to the first term only (CSL preconditioner), where only one V-cycle is required to obtain an (approximate) preconditioner inverse. These observations are comparable to the 1D results from Section \ref{sec:multires}.

In Table \ref{tab:EXmultigrid_2d_2} results for solving the same system using $EX(m)$-preconditioned GMRES are shown. Note that contrary to the BiCGStab results in Table \ref{tab:EXmultigrid_2d}, the $EX(m)$ preconditioner does appear to be cost-efficient for the 2D Helmholtz problem with wavenumber $k^2 = 2$e$+4$ for values of $m > 1$. Indeed, in this case the optimal preconditioner with respect to CPU time is the second-order polynomial $EX(2)$, which reduces the number of outer Krylov iterations to $191$ and minimizes the CPU time to $477.3$ seconds, compared to 540.8 seconds for $EX(1)$-GMRES. The main cause for this phenomenon is the relatively high per-iteration computational cost of the GMRES algorithm, which is caused by the orthogonalization procedure with respect to the Krylov subspace basis vectors. This cost is especially pronounced for larger iteration numbers. Hence, the good approximation properties of the $EX(m)$ preconditioner for higher values of $m$ may prove useful when the Krylov iteration cost is non-marginal compared to the cost of applying the preconditioner. The idea of polynomial preconditioners for GMRES has recently been proposed in the literature, see e.g.~\cite{liu2015polynomial}.

\begin{table}[t]
\centering {\small
\begin{tabular}{ C{0.8cm} | C{1.4cm} | C{2.0cm} | C{1.4cm} | C{2.0cm} |}
\cline{2-5}
& \multicolumn{2}{c| }{$n_x \times n_y = 128 \times 128$} & \multicolumn{2}{ c| }{$n_x \times n_y = 256 \times 256$} \\
& \multicolumn{2}{c| }{$k^2 = 5$e$+3$} & \multicolumn{2}{ c| }{$k^2 = 2$e$+4$} \\
\hline
\multicolumn{1}{|c||}{$m$} & $p(m)$ & CPU time & $p(m)$ & CPU time\\
\hline
\multicolumn{1}{|c||}{1} & 67 & 19.0 s. & 233 & 540.8 s. \\
\multicolumn{1}{|c||}{2} & 50 & 22.9 s. & 191 & 477.3 s. \\
\multicolumn{1}{|c||}{3} & 41 & 26.8 s. & 175 & 497.5 s. \\
\multicolumn{1}{|c||}{4} & 37 & 31.8 s. & 168 & 547.8 s. \\
\multicolumn{1}{|c||}{5} & 34 & 35.9 s. & 165 & 611.3 s. \\
\hline
& \multicolumn{4}{c| }{EX(m)-preconditioned GMRES} \\ 
\cline{2-5}
\end{tabular}
}
\vspace{0.3cm}
\caption{Performance of $EX(m)$-GMRES for different values of $m$ on the discretized 2D Helmholtz model problem \eqref{eq:model_problem_2d}. V(1,1)-cycle approximate preconditioner inversion. Column 1 \& 3: number of $EX(m)$-GMRES iterations $p(m)$. Col.~2 \& 4: total CPU time until convergence (in seconds).}
\label{tab:EXmultigrid_2d_2}
\end{table}

\section{Conclusions} 
\label{sec:conclusions}

In this work we have proposed a theoretical framework that generalizes the classic shifted Laplacian preconditioner by introducing the class of polynomial expansion preconditioners $EX(m)$. This concept extends the one-term CSL preconditioner to an $m$-term Taylor polynomial in the inverse complex shifted Laplace operator.  The outer iteration for solving the preconditioned system used in this work is a traditional Krylov iteration such as BiCGStab or GMRES.

Key properties of the $EX(m)$ preconditioner class are its structure as a finite $m$-term series of powers of CSL inverses (Neumann series), and its resulting asymptotic exactness, meaning $EX(m)$ approaches $A^{-1}$ in the limit for $m$ going to infinity. The polynomial structure of the $EX(m)$ preconditioner makes it easy to compute an iterative approximation to this polynomial, using e.g.\ $m$ multigrid V(1,1)-cycles (one for each term), which are guaranteed to converge given that the complex shift is chosen to be sufficiently large.

The preconditioning efficiency of the $EX(m)$ preconditioner is validated using a classic eigenvalue analysis. The addition of extra terms in the preconditioning polynomial clusters the spectrum around $1$, which reduces the condition number of the preconditioned Helmholtz operator and suggests a significant reduction in the number of outer Krylov iterations for large $m$. 

Numerical results on 1D and 2D Helmholtz benchmark problems support the theoretical results. The number of outer Krylov iterations is reduced significantly by the higher degree expansion preconditioners. Unfortunately, the computational cost of applying the $EX(m)$ preconditioner is directly proportionate to the number of terms $m$, since an extra (approximate) CSL inversion is required for each additional term in the polynomial. The use of a large number of series terms is thus not necessarily guaranteed to result in a more cost-efficient preconditioner. 

Following the numerical results of the 1D and 2D experiments, these conclusions are expected to be directly generalizable to higher spatial dimensions. Moreover, the constant wavenumber experiments performed in this paper are extensible to Helmholtz problems with heterogeneous and/or discontinuous wavenumbers, provided that the complex shift variable in the $EX(m)$ polynomial is large enough to allow for a stable numerical solution of the shifted Laplace systems for all wavenumber regimes occuring in the problem.
Hence, from a practical preconditioning interest, the simple one-term shifted Laplace preconditioner appears to be the optimal member of the $EX(m)$ class for many applications and problem configurations. 

Furthermore, two generalizations to the class of expansion preconditioners were presented and analyzed. These generalizations primarily prove to be insightful from a theoretical point of view. It is shown that the Taylor expansion preconditioner can be substituted by an optimal $m$-term polynomial which is theoretically cost-efficient. However, in practical applications the reduction of the number of outer Krylov iterations due to $EX(m)$ preconditioning often does not pay off to the cost of the extra approximate multigrid inversions.

For systems in which the cost of applying the outer Krylov step becomes significant relative to the cost of the preconditioner application, the use of multiple terms in the $EX(m)$ expansion could result in a more cost-efficient solver. Possible scenarios for this include the application of non-restarted GMRES as the outer Krylov solver, the use of higher order discretization schemes for the original Laplace operator (while maintaining second order discretization for the preconditioner), etc. 

Additionally, the treatment of extremely large-scale HPC systems on massively parallel hardware may warrant the need for higher-order polynomial preconditioners. Since Krylov methods are typically communication (or bandwidth) bound instead of compute bound in this context, polynomial preconditioning directly reduces the number of communication bottlenecks (dot-products) by reducing the number of Krylov iterations, while simultaneously improving the arithmetic intensity of the solver. A detailed analysis of these individual scenarios is however well beyond the scope of this text, and is left for future work.

The generalization to the shifted Laplace preconditioner for Helmholtz problems proposed in this work is particularly valuable from a theoretical viewpoint, providing fundamental insights into the concept of shifted Laplace preconditioning by situating the classic complex shifted Laplace operator in a broader theoretical context, and proving the classic CSL preconditioner to be the most cost-efficient member of the $EX(m)$ preconditioner class for the most common practical problems.

\section*{Acknowledgments}

This research in this chapter was funded by the Research Council of the University of Antwerp under BOF grant number 41-FA070300-FFB5342. The scientific responsibility rests with its authors.
The authors would like to cordially thank dr.~Bram Reps for fruitful discussions on the subject.

\nocite{*}
{\footnotesize
\bibliographystyle{plain}
\bibliography{refs}
}

\end{document}